\documentclass[11pt]{article}
\usepackage[T1]{fontenc}
\usepackage{natbib}
\usepackage{amsmath,amsfonts,amssymb,amsthm}
\usepackage{a4wide}
\usepackage{color}
\usepackage{epsfig,psfrag}

\newtheorem{dfn}{Definition}[section]
\newtheorem{thm}[dfn]{Theorem}
\newtheorem{lem}[dfn]{Lemma}
\newtheorem{cor}[dfn]{Corollary}
\newtheorem{rem}[dfn]{Remark}
\newtheorem{prop}[dfn]{Proposition}

\newenvironment{dem}{\vskip 2mm\noindent {\it Proof} :}
                    {\hfill $\square$ \vskip 2mm \noindent}
 
\DeclareMathOperator{\cyl}{cyl}
\DeclareMathOperator{\hyp}{hyp}
\DeclareMathOperator{\card}{card}

\DeclareMathOperator{\ad}{ad}

\DeclareMathOperator{\argmin}{argmin}

\DeclareMathOperator{\lleft}{left}
\DeclareMathOperator{\rright}{right}

\newcommand{\eps}{\varepsilon}

\def\PP{\mathbb{P}}
\def\RR{\mathbb{R}}
\def\EE{\mathbb{E}}
\def\NN{\mathbb{N}}
\def\ZZ{\mathbb{Z}}
\def\E{\mathcal{E}}
\def\D{\mathcal{D}}

\def\N{\mathcal{N}}
\def\I{\mathcal{I}}
\def\J{\mathcal{J}}
\def\M{\mathcal{M}}
\def\G{\mathcal{G}}
\def\K{\mathcal{K}_{\theta,h}}
\def\wK{\widetilde{\mathcal{K}}_{\theta,h}}
\def\wt{{\widetilde{\theta}}}
\def\F{\mathcal{F}}

\def\LL{\mathbb{L}}
\def\va{\vec{v}(\theta)}
\def\vb{\vec{v}(\widetilde{\theta})}
\def\vc{\vec{v}^{\bot}(\theta)}
\def\vd{\vec{v}^{\bot}(\widetilde{\theta})}

\def\vg{\vec{v}^{\bot}(\widetilde{\theta}_b)}
\def\vh{\vec{v}^{\bot}(\widetilde{\theta}_a)}

\begin{document}

\thispagestyle{empty}


\title{Lower large deviations for the
  maximal flow through tilted cylinders in two-dimensional first passage
  percolation}

\date{}
\author{}
\maketitle

\begin{center}
\vskip-1cm {\Large Rapha\"el Rossignol}\\
{\it Universit\'e Paris Sud, Laboratoire de Math\'ematiques, b\^atiment 425,
91405 Orsay Cedex, France}\\
{\it E-mail:} raphael.rossignol@math.u-psud.fr\\
and\\
\vskip1cm {\Large Marie Th\'eret}\\
{\it \'Ecole Normale Sup\'erieure, D\'epartement Math\'ematiques et
Applications, 45 rue d'Ulm, 75230 Paris Cedex 05, France}\\
{\it E-mail:} marie.theret@ens.fr
\end{center}



\begin{abstract}
Equip the edges of the lattice $\ZZ^2$ with i.i.d. random
capacities. A law of large numbers is known for
the maximal flow crossing a rectangle in $\RR^2$ when the side
lengths of the rectangle go to infinity. We prove that the lower large
deviations are of surface order, and we prove the corresponding large
deviation principle from below. This extends and improves
previous large deviations results of \cite{GrimmettKesten84} obtained
for boxes of particular orientation.
\end{abstract}

\noindent
{\it AMS 2000 subject classifications:} Primary 60K35; secondary 82B43.

\noindent
{\it Keywords :} First passage percolation, maximal flow, large deviation
principle.

\section{Introduction}

Imagine each edge of $\ZZ^2$ is a microscopic pipe some fluid can go
through. To each edge $e$, we attach a random capacity $t(e)$ and
suppose that all these random variables are independent and
identically distributed with common distribution function $F$ on $\RR^+$. Now,
we take a large rectangle $R$ in $\RR^2$, decide that one side is the
``left side'', and accordingly name the other sides the right side, the
top and the bottom of $R$. We are concerned with the \emph{maximal
  flow rate} that
can cross $R$ from the top to the bottom, while never exceeding the
capacities of the edges (see section
\ref{chapitre6subsec:squarelattice} for a formal definition). Informally, we ask:  how
does the maximal flow between the top and the bottom behave when $R$
 gets larger and larger ? This question was first considered in
 \cite{GrimmettKesten84}, where a law of large numbers and large deviation estimates where proved, but
only for ``straight'' rectangles $R$, i.e. with sides parallel to the
coordinate axes. Let us mention that lower large deviations are of surface
order, i.e. of the order of the width of $R$ whereas upper large
deviations are of volume order, i.e. of the order of the area of $R$. In a previous work, \cite{RossignolTheret09}, the authors
extended the law of large numbers of \cite{GrimmettKesten84} to
rectangles $R$ with arbitrary orientation. The purpose of this article is to give
a corresponding lower large deviation principle. 

We shall precise now our contribution. Let us notice that the problem of maximal flow has been
studied in the more general case of the lattice $\ZZ^d$,  $R$ being
then some box in $\RR^d$, first through the work of \cite{Chayes} and
\cite{Kesten:flows}, and then notably with \cite{Zhang},
\cite{Zhang07}, \cite{TheretUpper}, \cite{Theret:small},
\cite{RossignolTheret08b}. In any case, it is much simpler to study
a subadditive version of the maximal flow, which we shall call $\tau$
in section \ref{chapitre6subsec:squarelattice} below. Then, when $R$ is
straight (i.e. has faces parallel to coordinate hyperplanes), one can pass from $\tau$ to the maximal flow thanks to
symmetry considerations. When  $d\geq 3$,
results concerning $R$ such as large deviations estimates or laws of
large numbers and that do not suppose that $R$ is straight (or even
that it is a box) are known from the recent works \cite{CerfTheret09geo},
\cite{CerfTheret09inf} and \cite{CerfTheret09sup}, but they require a
lot of geometric work, need strong moment hypotheses, have not yet
provided large deviation principles and in any case
are much more involved than what is expected to be necessary in two
dimensions to treat the simple case of rectangles. Indeed, in two dimensions, duality considerations are of
great help to prove a law of large numbers and a lower large deviation
principle. The aim of  \cite{RossignolTheret09} was thus to
derive in a simple way the law of large numbers for the maximal flow
from the top to the bottom through a rectangle $R$. The
aim of the present article is to use the same constructions to derive
in a simple way the lower large deviation principle of the maximal
flow. A similar, simple proof in dimension $d\geq 3$ remains to be
found. 


The precise definitions and notations, the relevant background as well
as the main results are presented in section
\ref{chapitre6sec:notations}. Then, the main construction and the
lower large deviation estimates are the purpose of section \ref{chapitre6subsec:deviations}, while
section \ref{chapitre6subsec:LDP} contains the proof of the large deviation principle itself.


\section{Notations,  background and main results}
\label{chapitre6sec:notations}

The most important notations are gathered in sections
\ref{chapitre6subsec:maxflow} to \ref{chapitre6subsec:duality}, the
relevant background is described in section \ref{secbackground} while our
main results are stated in section \ref{secmainres}. We discuss in section
\ref{remcond} the different hypotheses appearing in our results.

\subsection{Maximal flow on a graph}
\label{chapitre6subsec:maxflow}
First, let us define the notion of a flow on a finite
unoriented graph $G=(V,\E)$ with set of vertices $V$ and set of edges
$\E$.  Let
$t=(t(e))_{e\in \E}$ be a collection of non-negative real numbers,
which are called capacities. It means that $t(e)$ is the maximal amount of
    fluid that can go through the edge $e$ per unit of time. To each edge $e$, one may associate two oriented edges, and we
shall denote by $\smash{\overrightarrow{\E}}$ the set of all these oriented edges. Let $A$ and $Z$ be
two finite, disjoint, non-empty sets of vertices of $G$: $A$ denotes the source of the
network, and $Z$ the sink. A function $\theta$
on $\smash{\overrightarrow{\E}}$ is called a flow from $A$ to $Z$ with strength
$\|\theta\|$ and capacities $t$ if it is
antisymmetric, i.e.
$\theta_{\overrightarrow{xy}}=-\theta_{\overrightarrow{yx}}$, if it satisfies the
node law at each vertex $x$ of $V\smallsetminus (A\cup Z)$:
$$\sum_{y\sim x}\theta_{\overrightarrow{xy}}=0\;,$$
where $y\sim x$ means that $y$ and $x$ are neighbours on $G$,
if it satisfies the capacity constraints:
$$\forall e\in \E,\;|\theta(e)|\leq t(e)\;,$$
and if the ``flow in'' at $A$ and the ``flow out'' at $Z$ equal $\|\theta\|$:
$$\|\theta\|=\sum_{a\in A}\sum_{\substack{y\sim a\\ y\not \in
    A}}\theta(\overrightarrow{ay})=\sum_{z\in Z}\sum_{\substack{y\sim
    z\\ y\not \in Z}}\theta(\overrightarrow{yz})\;.$$
The \emph{maximal flow from $A$ to $Z$}, denoted by $\phi_t(G,A,Z)$, is defined as
the maximum strength of all flows from $A$ to $Z$ with capacities
$t$. We shall in general omit the subscript $t$ when it is understood
from the context. The \emph{max-flow min-cut theorem} (see \cite{Bollobas} for instance)
asserts that the maximal flow from $A$ to $Z$ equals the minimal
capacity of a cut between $A$ and $Z$. Precisely, let us say that
$E\subset\E$ is a cut between $A$ and $Z$ in $G$ if every path from
$A$ to $Z$ borrows at least one edge of $E$. Define $V(E)=\sum_{e\in
  E}t(e)$ to be the capacity of a cut $E$. Then,
\begin{equation}
\label{chapitre6eq:maxflowmincut}
\phi_t(G,A,Z)=\min\{V(E)\mbox{ s.t. }E\mbox{ is a cut between
}A\mbox{ and }Z \mbox{ in } G\}\;.
\end{equation}


\subsection{On the square lattice}
\label{chapitre6subsec:squarelattice}
We shall always consider $G$ as a piece of $\ZZ^2$. More precisely, we consider the graph $\LL=(\mathbb{Z}^{2},
\mathbb E ^{2})$ having for vertices $\mathbb Z ^{2}$ and for edges
$\mathbb E ^{2}$, the set of pairs of nearest neighbours for the standard
$L^{1}$ norm. The notation $\langle x,y\rangle$ corresponds to the edge with endpoints
$x$ and $y$. To each edge $e$ in $\mathbb{E}^{2}$ we associate a random
variable $t(e)$ with values in $\mathbb{R}^{+}$. \emph{We suppose that the family
$(t(e), e \in \mathbb{E}^{2})$ is independent and identically distributed,
with a common distribution function $F$}. More
formally, we take the product
measure $\mathbb {P}=F^{\otimes \Omega}$ on $\Omega= \prod_{e\in \mathbb{E}^{2}} [0, \infty[$,
    and we write its expectation $\mathbb{E}$.  If $G$ is
    a subgraph of $\LL$, and $A$ and $Z$ are two subsets of vertices of
    $G$, \emph{we shall denote by $\phi(G,A,Z)$ the maximal flow in $G$ from $A$ to
    $Z$}, where $G$ is equipped with capacities $t$. When $B$ is a
    subset of $\RR^2$, and $A$ and $Z$ are subsets of $\ZZ^2\cap B$,
    we shall denote by $\phi(B,A,Z)$ again the maximal flow
    $\phi(G,A,Z)$ where $G$ is the induced subgraph of $\ZZ^2$ with set of vertices
    $\ZZ^2\cap B$.


We denote by
$\overrightarrow{e}_1$ (resp. $\overrightarrow{e}_2$) the vector $(1,0)\in\RR^2$ (resp. $(0,1)$). Let $A$ be a non-empty
line segment in $\RR^2$. We shall denote by $l(A)$ its (Euclidean)
length. All line segments will be
supposed to be closed in $\mathbb{R}^2$. We denote by
$\va$ the vector of unit Euclidean norm orthogonal to $\hyp (A)$,
the hyperplane spanned by $A$, and such that there is
$\theta\in[0,\pi[$ such that $\va=(\cos \theta,\sin\theta)$. Define $\vc=( \sin\theta,-\cos\theta)$ and
denote by $a$ and $b$ the end-points of $A$ such that
$(b-a).\vc >0$. For $h$ a positive real number, \emph{we denote by $\cyl(A,h)$ the cylinder of basis $A$ and height $2h$},
i.e., the set 
$$ \cyl (A,h) \,=\, \{x+t \va \,|\, x\in A \,,\,  t\in
[-h,h]    \}\,.$$
We define also \emph{the $r$-neighbourhood $\mathcal{V} (H,r)$ of a subset $H$ of $\mathbb{R}^d$} as
$$\mathcal{V}(H,r) \,=\, \{ x \in \mathbb{R}^d \,|\, d(x,H)<r\}\,,$$
where the distance is the Euclidean one ($d(x,H) = \inf \{\|x-y\|_2
\,|\, y\in H \}$).

\begin{figure}[!ht]
\centering
\input{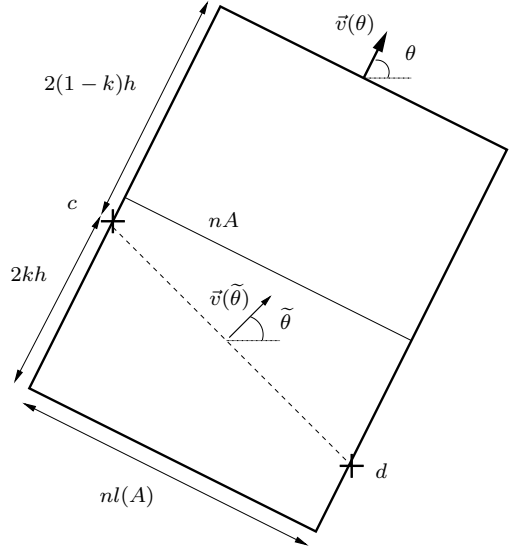}
\caption{An admissible boundary condition $(k,\wt)$.}
\label{chapitre6fig:notations}
\end{figure}

Now, we define $D(A,h)$ the set of \emph{admissible boundary conditions} on
$\cyl(A,h)$ (see Figure \ref{chapitre6fig:notations}): 
$$D(A,h)=\left\{(k,\tilde \theta) \,|\, k\in [0,1]\mbox{ and
  }\tilde\theta\in
  \left[\theta-\arctan\left(\frac{2hk}{l(A)}\right),\theta+\arctan\left(\frac{2h(1-k)}{l(A)}\right)\right]\right\}\;.$$
The meaning of an element $\kappa=(k,\tilde \theta)$ of $D(A,h)$ is the
following. We define
$$\vb \,=\, (\cos \widetilde{\theta}, \sin
\widetilde{\theta}) \qquad \textrm{and} \qquad \vd
\,=\, (\sin \widetilde{\theta}, -\cos \widetilde{\theta})\,. $$
In $\cyl(nA,h(n))$, we may define two points $c$ and $d$
such that $c$ is ``at height $2kh$ on the left side of $\cyl(A,h)$'',
and $d$ is ``on the right side of $\cyl(A,h)$''  by 
$$c=a+(2k-1)h\va \,,\quad
(d-c)\mbox{ is orthogonal to }\vb \quad\mbox{ and
}d\mbox{ satisfies }\vec{cd}\cdot \vd >0 \;.$$
Then we see that $D(A,h)$ is exactly the set of parameters so that $c$ and $d$
remain ``on the sides of $\cyl(A,h)$''. 

We define also $\D(A,h)$, the
set of angles $\tilde\theta$ such that there is an admissible boundary
condition with angle $\tilde\theta$:
$$\D(A,h)=\left[\theta-\arctan\left(\frac{2h}{l(A)}\right),\theta+\arctan\left(\frac{2h}{l(A)}\right)\right]\;.$$
It will be useful to define
 the \emph{left side (resp. right side) of $\cyl(A,h)$}: let $\lleft(A)$
 (resp. $\rright(A)$) be the set of vertices in $\cyl(A,h)\cap\ZZ^2$
such that there exists $y\notin \cyl(A,h)$, $\langle x,y\rangle \in
\mathbb{E}^d$ and $[ x,y [$, the segment that includes $x$ and excludes
$y$, intersects $a+[-h,h].\va $ (resp. $b+[-h,h].\va $).

Now, the set $\cyl(A,h) \smallsetminus (c+\RR(d-c))$ has two connected
components, which we denote by $\mathcal{C}_1(A,h,k,\tilde\theta)$ and
$\mathcal{C}_2(A,h,k,\tilde\theta)$. For $i=1,2$, let $A_i^{h,k,\tilde\theta}$ be
the set of the points in $\mathcal{C}_i(A,h,k,\tilde\theta) \cap \mathbb{Z}^2$ which have
a nearest neighbour in $\mathbb{Z}^2 \smallsetminus \cyl(A,h)$:
$$A_i^{h,k,\tilde\theta}\,=\,\{x\in \mathcal{C}_i(A,h,k,\tilde\theta) \cap
\mathbb{Z}^2 \,|\, \exists y \in \mathbb{Z}^2 \smallsetminus \cyl(A,h) \,,\,
\|x-y\|_{1} =1 \}\,.$$
We define \emph{the flow in $\cyl(A,h)$ constrained by the boundary
condition $\kappa=(k,\tilde\theta)$} as:
$$\phi^\kappa(A,h):=\phi(\cyl(A,h),A_1^{h,k,\tilde\theta},A_2^{h,k,\tilde\theta})\;.$$
A special role is played by the condition  $\kappa = (1/2, \theta)$,
and we shall denote: 
$$\tau(A,h) = \tau (\cyl(A,h), \va ):=\phi^{(1/2, \theta)}(A,h)\;.$$
Let $T(A,h)$ (respectively $B(A,h)$) be the top
(respectively the bottom) of $\cyl(A,h)$, i.e.,
$$ T(A,h) \,=\, \{ x\in \cyl(A,h) \,|\, \exists y\notin \cyl(A,h)\,,\,\,
\langle x,y\rangle \in \mathbb{E}^d \mbox{ and }\langle x,y\rangle
\mbox{ intersects } A+h\va  \}  $$
and
$$  B(A,h) \,=\, \{ x\in \cyl(A,h) \,|\, \exists y\notin \cyl(A,h)\,,\,\,
\langle x,y\rangle \in \mathbb{E}^d \mbox{ and }\langle x,y\rangle
\mbox{ intersects } A-h\va \} \,.$$
We shall denote the flow in $\cyl(A,h)$ from the top to the bottom as:
$$\phi(A,h)= \phi (\cyl(A,h), \va)=\phi(\cyl(A,h),T(A,h),B(A,h))\;.$$


\subsection{Duality}
\label{chapitre6subsec:duality}

The main reason why dimension 2 is easier to deal with than dimension
$d\geq 3$ is duality. Planar duality implies that there are only
$O(h^2)$ admissible boundary conditions on $\cyl(A,h)$, as shows the
following lemma taken from \cite{RossignolTheret09}.
\begin{lem}
\label{chapitre6lem:duality}
Let $A$ be any line segment in $\RR^2$ and $h$ a positive real
number. Then,
\begin{equation*}
\phi(A,h)= \min_{\kappa\in D(A,h)}\phi^\kappa(A,h) \;.
\end{equation*}
\end{lem}
Notice that the condition $\kappa$ belongs to the non-countable set
$D(A,h)$, but the graph is discrete so $\phi^\kappa(A,h)$ takes only a
finite number of values when $\kappa\in D(A,h)$. Precisely, there is a
finite subset $\tilde D(A,h)$ of $D(A,h)$, such that:
\begin{equation}
\label{chapitre6eq:Oh2}
\card(\tilde D(A,h))\leq C_4 h^2\;,
\end{equation}
for some universal constant $C_4$, and:
\begin{equation*}
\phi(A,h)= \min_{\kappa\in \tilde D(A,h)}\phi^\kappa(A,h) \;.
\end{equation*}


\subsection{Background}
\label{secbackground}
We gather in this section known results concerning the behaviour of the variables $\tau(nA,h(n))$ and $\phi(nA,h(n))$ when $n$ and $h(n)$ go to
infinity. They are of two types. We present first the law of large numbers
satisfied by both variables. Then we describe the large deviation principle
from below proved for $\tau(nA,h(n))$. The purpose of this article is to
extend the study of lower large deviations to the variable $\phi(nA,h(n))$.

We gather the main hypotheses that we shall do on $F$
and on the height $h$. Notice that $\mathbf{ (F5)}\Rightarrow\mathbf{ (F4)}\Rightarrow\mathbf{
  (F3)}\Rightarrow\mathbf{ (F2)}$ and $\mathbf{
  (H3)}\Rightarrow\mathbf{ (H1)}$.
{\begin{center}
\renewcommand{\arraystretch}{1.5}
\begin{tabular}{|l|l|}
\hline
\multicolumn{1}{|c|}{Hypotheses on $F$} &
\multicolumn{1}{|c|}{Hypotheses on $h$}
\\
\hline
$\mathbf{(F1)}\quad F(0)<1-p_c(d)$&    $\mathbf{ (H1)}\quad \lim_{n
  \rightarrow \infty}h(n)=+\infty$ \\
$\mathbf{ (F2)}\quad \int_0^\infty
x\;dF(x)<\infty$ & $\mathbf{ (H2)}\quad \lim_{n \rightarrow
  \infty}\frac{\log h(n)}{n}=0$ \\ 
$\mathbf{ (F3)}\quad \int_0^\infty x^{2}\;dF(x)<\infty$ & $\mathbf{
  (H3)}\quad \lim_{n \rightarrow \infty}\frac{h(n)}{n}=+\infty$  \\
$ \mathbf{ (F4)}\quad \exists \gamma>0 ,\;  \int_0^\infty e^{\gamma
  x}\;dF(x)<\infty $ & $ \mathbf{(H4)}\quad \exists \alpha \in
\left[0,\frac{\pi}{2}\right] ,\;\lim_{n\rightarrow \infty}
\frac{2h(n)}{nl(A)}=\tan \alpha $   \\
$ \mathbf{ (F5)}\quad \forall \gamma>0 ,\;  \int_0^\infty e^{\gamma
  x}\;dF(x)<\infty $ & \\
\hline 
\end{tabular}
\end{center}}
\vspace*{2mm}
Using a subadditive argument and deviation inequalities, Rossignol
and Th\'eret have proved in \cite{RossignolTheret08b} that $\tau(nA,
h(n))$ satisfies a law of large numbers:
\begin{thm}
\label{chapitre6thm:LGNtau}
We suppose that $\mathbf{(F2)}$ holds.
For every unit vector $\va =(\cos\theta,\sin\theta)$, there exists a
constant $\nu_\theta$ depending on $F$ and $\theta$, such that for
every non-empty line-segment $A$ orthogonal to $\va$ and of Euclidean length $l(A)$, for
every height function $h: \NN \rightarrow \RR^+$ satisfying $\mathbf{(H1)}$, we have
$$ \lim_{n\rightarrow \infty} \frac{\tau(nA, h(n))}{n l(A)}
\,=\, \nu_\theta \qquad \textrm{in } L^1 \,. $$
Moreover, if the origin of the graph belongs to $A$, or if $\mathbf{(F3)}$
holds, then
$$ \lim_{n\rightarrow \infty} \frac{\tau(nA, h(n))}{nl(A)}
\,=\, \nu_\theta \qquad \textrm{a.s.} $$
\end{thm}
This law of large numbers holds in fact for every dimension $d\geq 2$, and
the limit $\nu$ depends also on the dimension. Let
us remark that (in dimension two) $\nu_\theta$ is equal to $\mu(\vc )$, where
$\mu(.)$ is the time-constant function of first passage percolation as
defined in \cite{Kesten:StFlour}, (3.10) p.~158. This equality follows
from the duality considerations of section \ref{chapitre6subsec:duality} and
standard first passage percolation techniques (see also
Theorem 5.1 in \cite{GrimmettKesten84}) that relate cylinder passage times to
unrestricted passage times (as in \cite{HammersleyWelsh},
Theorem~4.3.7 for instance). Boivin has also proved a very similar law of large
numbers (see Theorem 6.1 in \cite{Boivin}). Notice that for the
definition of $\mu(.)$, Kesten requires only the existence of the first
moment of the law $F$ in the proof from \cite{Kesten:StFlour}, and it
can also be defined under the weaker condition $\int_0^{\infty}(1
-F(x))^4\;dx<\infty$. 

One consequence of this equality between $\nu$ and $\mu$ is that
$\theta\mapsto \nu_\theta$ is either constant equal to zero, or always
non-zero. In fact the following property holds (cf. \cite{Kesten:StFlour},
Theorem 6.1 and Remark 6.2 p.~218):
\begin{prop}
\label{propnu}
We suppose that $\mathbf{(F2)}$ holds. Then $\nu_\theta$
is well defined for all $\theta$, and we have
$$ \nu_\theta \,>\,0 \iff \mathbf{(F1)}\,.  $$
\end{prop}

We recall that for all $n\in \NN$, we have defined
$$\D(nA,h(n))=\left[\theta-\arctan\left(\frac{2h(n)}{nl(A)}\right),\theta+\arctan\left(\frac{2h(n)}{nl(A)}\right)\right]\;.$$
Extending the law of large numbers proved by \cite{GrimmettKesten84} for
$\phi(nA,h(n))$ in boxes of particular orientation, the authors proved the
following result (see Theorem 2.8 and Corollary 2.10 in
\cite{RossignolTheret09}), in the same spirit as the result of \cite{Garet2}:
\begin{thm}
\label{chapitre6thm:lgn}
Let $A$ be a non-empty line-segment in $\RR^2$, with Euclidean length
$l(A)$. Let $\theta\in[0,\pi[$ be such that $(\cos\theta,\sin\theta)$
is orthogonal to $A$ and $h:\NN \rightarrow \RR^+$ satisfying
$\mathbf{(H1)}$ and $\mathbf{(H2)}$. Define:
$$\overline{\D}=\limsup_{n\rightarrow
  \infty}\D(nA,h(n))=\bigcap_{N\geq 1}\bigcup_{n\geq
  N}\D(nA,h(n))\;,$$
and
$$\underline{\D}=\liminf_{n\rightarrow \infty}\D(nA,h(n))=\bigcup_{N\geq 1}\bigcap_{n\geq N}\D(nA,h(n))\;.$$
Suppose that $\mathbf{(F2)}$ holds. Then,
\begin{equation}
\label{eq:liminfmoyenne}
\liminf_{n\rightarrow \infty}\frac{\EE[\phi(nA,h(n))]}{nl(A)}=\inf \left\{
  \frac{\nu_{\widetilde{\theta}}}{\cos  (\widetilde{\theta} - \theta)}
  \,|\, \widetilde{\theta} \in \overline{\D} \right\}
\end{equation}
and
\begin{equation}
\label{eq:limsupmoyenne}\limsup_{n\rightarrow\infty}\frac{\EE[\phi(nA,h(n))]}{nl(A)}= \inf\left\{
  \frac{\nu_{\widetilde{\theta}}}{\cos (\widetilde{\theta} - \theta)} \,|\,
  \widetilde{\theta} \in \underline{\D}  \right\} \,.
\end{equation}
Moreover, if the origin of the graph $0$ is the middle of $A$, or if
$\mathbf{(F3)}$ holds, then
$$\liminf_{n\rightarrow \infty}\frac{\phi(nA,h(n))}{nl(A)}=\inf \left\{
  \frac{\nu_{\widetilde{\theta}}}{\cos  (\widetilde{\theta} - \theta)}
  \,|\, \widetilde{\theta} \in \overline{\D} \right\} \qquad a.s.$$
and
$$\limsup_{n\rightarrow\infty}\frac{\phi(nA,h(n))}{nl(A)}= \inf\left\{
  \frac{\nu_{\widetilde{\theta}}}{\cos (\widetilde{\theta} - \theta)} \,|\,
  \widetilde{\theta} \in \underline{\D}  \right\} \qquad a.s.$$
\end{thm}
\begin{cor}
\label{chapitre6corollaire}
We suppose that conditions $\mathbf{(H1)}$, $\mathbf{(H2)}$ and
$\mathbf{(H3)}$ on $h$ are satisfied. We suppose also that $\mathbf{(F2)}$
holds. Then we have
$$\lim_{n\rightarrow \infty}\frac{\phi(nA,h(n))}{nl(A)}=\inf\left\{
  \frac{\nu_{\widetilde{\theta}}}{\cos (\widetilde{\theta} - \theta)} \,|\,
  \widetilde{\theta} \in [\theta-\alpha,\theta+\alpha] \right\} \qquad
\textrm{in } L^1\,.$$
Moreover, if $0$ is the middle of $A$, or if
$\mathbf{(F3)}$ holds, then
$$\lim_{n\rightarrow \infty}\frac{\phi(nA,h(n))}{nl(A)}=\inf\left\{
  \frac{\nu_{\widetilde{\theta}}}{\cos (\widetilde{\theta} - \theta)} \,|\,
  \widetilde{\theta} \in [\theta-\alpha,\theta+\alpha] \right\} \qquad
\textrm{a.s.}$$
\end{cor}

Concerning the lower large deviations of $\tau(nA,h(n))$, Theorem 3.9 and
Lemma 5.1 in \cite{RossignolTheret08b} state that:
\begin{thm}
\label{propI}
For every non-empty line-segment $A$ in $\RR^2$, with Euclidean length
$l(A)$, for every height function $h: \NN \rightarrow \RR^+$ satisfying
$\mathbf{(H1)}$, for all $\lambda$ in $\RR^+$, the limit
$$ \I_{\theta} (\lambda) \,=\,\lim_{n\rightarrow \infty} \frac{-1}{nl(A)} \log \PP \left[ \tau (nA,
  h(n) ) \leq \left( \lambda - \frac{1}{\sqrt n}  \right) n l(A)  \right] $$
exists in $[0,+\infty]$ and depends only on $\theta\in[0,\pi[$ such that
$(\cos\theta,\sin\theta)$ is orthogonal to $A$, and not on $h$ nor $A$
itself. Moreover, if the hypotheses $\mathbf{(F1)}$ and $\mathbf{(F2)}$ are
satisfied, the function $\I_{\theta}$ has the following properties: it is
convex on $\RR^+$, infinite on $[0,\delta (|\cos \theta| + |\sin
\theta|)[$, where $\delta = \inf \{\lambda \,|\, \PP [t(e) \leq \lambda]>0
\}$, finite on $]\delta(|\cos \theta| + |\sin \theta|), +\infty[ $, equal
to $0$ on $[\nu_{\theta}, +\infty[$, and if $\nu_{\theta} > \delta(|\cos
\theta| + |\sin \theta|) $ it is continuous and strictly decreasing on $]\delta(|\cos \theta| + |\sin \theta|),
\nu_{\theta}]$ and positive on $]\delta(|\cos \theta| + |\sin \theta|),
\nu_{\theta}[$.
\end{thm}

For simplicity of notations, we define $\I_{\wt}=+\infty $ on $\RR^-_*$,
and for all $a\geq 0$,
$$\I_{\wt} (a^+) \,=\, \lim_{\varepsilon
  \rightarrow 0, \,\varepsilon >0} \I_{\wt} (a+\varepsilon) \qquad
\mbox{and} \qquad
\I_{\wt}(a^-) \,=\, \lim_{\varepsilon \rightarrow 0,\,\varepsilon >0} \I_{\wt}
(a-\varepsilon)\,.$$
We denote by $\J_{\theta}$ the function defined on $\RR^+$ by
$$ \J_{\theta}(\lambda) \,=\, \left\{ \begin{array}{ll} \I_{\theta}
    (\lambda ^+) & if\,\, \lambda \leq \nu_{\theta}\,,\\ +\infty & if \,\,
    \lambda > \nu_{\theta} \,.\end{array} \right. $$
The following large deviation principle has also been proved in
\cite{RossignolTheret08b}, Theorem 3.10:
\begin{thm}
\label{thmPGDtau}
For every non-empty line-segment $A$ in $\RR^2$, with Euclidean length
$l(A)$, for every height function $h: \NN \rightarrow \RR^+$ satisfying
$\mathbf{(H1)}$, if $\mathbf{(F1)}$ and $\mathbf{(F5)}$ hold, then the
sequence
$$ \left( \frac{\tau(nA,h(n))}{nl(A)},n\in \NN \right) $$
satisfies a large deviation principle of speed $nl(A)$ with the good rate
function $\J_{\theta}$, where $\theta\in[0,\pi[$ is such that
$(\cos\theta,\sin\theta)$ is orthogonal to $A$.
\end{thm}

The same large deviation principle is also proved for $\phi(nA,h(n))$ if $\theta =0$ (see Theorem 3.17 in \cite{RossignolTheret08b},
condition $\mathbf{(F5)}$ is replaced by $\mathbf{(F4)}$ in comparison with
Theorem \ref{thmPGDtau}) or if $h$ satisfies $\lim_{n\rightarrow
  \infty}h(n)/n =0$ (see Corollary 3.14 in
\cite{RossignolTheret08b}). Theorems \ref{propI} and \ref{thmPGDtau} are valid in any dimension
$d\geq 2$. 
The difference of hypotheses
between theorems concerning the variable $\tau$ and theorems concerning the
variable $\phi$ will be discussed in Remark \ref{remcond}.


\subsection{Main results}
\label{secmainres}

As we have seen in Theorem \ref{chapitre6thm:lgn}, the existence of a limit
for $\phi(nA,h(n))$ is linked with the equality between different
infimum. The same holds for the large deviation principle, so we define two
additional hypotheses we will use:

{\begin{center}
\renewcommand{\arraystretch}{1.5}
\begin{tabular}{|l|}
\hline
\multicolumn{1}{|c|}{Hypotheses on $h$ and $F$} \\
\hline
$\mathbf{(FH1)}\quad \inf \left\{\frac{\nu_{\wt}}{\cos(\wt - \theta)} \,|\,
\wt \in \overline{\D} \right\} = \inf \left\{ \frac{\nu_{\wt}}{\cos(\wt -
  \theta)}\,|\, \wt \in \underline{\D} \right\} $ \\
$ \begin{array}{l} \hspace{-0.15cm} \mathbf{ (FH2)}\quad \forall \lambda \geq 0,\;  \inf
\left\{ \frac{1}{\cos(\wt - \theta)} \I_{\wt} (\lambda \cos(\wt - \theta)
  ^+) \,|\,\wt \in \ad(\overline{\D}) \right\} \\ \qquad \qquad \qquad \qquad \qquad  =
\inf \left\{ \frac{1}{\cos(\wt - \theta)} \I_{\wt} (\lambda \cos(\wt -
  \theta) ^+)  \,|\, \wt \in \ad(\underline{\D}) \right\} \end{array} $ \\
\hline 
\end{tabular}
\end{center}}
\vspace*{2mm}

Here and in the rest of the paper we denote the adherence of a set $S$ by $\ad(S)$. Notice that $\mathbf{(H4)}$ implies $\mathbf{(FH1)}$ and $\mathbf{(FH2)}$. We can now state our main results:
\begin{thm}[Lower Large Deviations]
\label{thmLDphi}
Let $A$ be a non-empty line-segment in $\RR^2$, $\theta \in [0,\pi[ $ such that
$(\cos \theta, \sin \theta)$ is orthogonal to $A$, and $h:\NN \rightarrow
\RR^+$satisfying conditions $\mathbf{(H1)}$ and $\mathbf{(H2)}$. If
$\mathbf{(F1)}$, $\mathbf{(F2)}$ and $\mathbf{(FH1)}$ hold, then there
exist constants $K_1(F,A,h,\eps)\geq 0$ and $K_2(F,\theta,h, \eps )>0$ such
that
$$ \PP\left( \frac{\phi(nA,h(n))}{nl(A)} \leq \eta_{\theta,h} -\eps \right)
\,\leq\, K_1 e^{-K_2 nl(A)}\,, $$
where
$$\eta_{\theta,h}\, =\, \inf_{\wt \in \overline{\D}} \frac{\nu_{\wt}}{ \cos
  (\wt - \theta)} \,=\, \lim_{n\rightarrow \infty}
\frac{\phi(nA,h(n))}{nl(A)} \qquad \textrm{in }L^1\,. $$
\end{thm}

\begin{thm}[Large Deviation Principle]
\label{thmPGDphi}
Let $A$ be a non-empty line-segment in $\RR^2$, $\theta \in [0,\pi[ $ such that
$(\cos \theta, \sin \theta)$ is orthogonal to $A$, and $h:\NN \rightarrow
\RR^+$satisfying conditions $\mathbf{(H1)}$ and $\mathbf{(H2)}$. If
$\mathbf{(F1)}$, $\mathbf{(F2)}$, $\mathbf{(FH1)}$, $\mathbf{(FH2)}$, and
either $\mathbf{(F4)}$ or $\mathbf{(H3)}$ hold, then the sequence
$$ \left( \frac{\phi(nA,h(n))}{nl(A)},n\in \NN \right) $$
satisfies a large deviation principle of speed $nl(A)$ with the good rate
function $\K: \RR^+ \rightarrow \RR^+ \cup \{+\infty\}$ defined by
$$ \K(\lambda) \,=\,\left\{ \begin{array}{ll} \inf \left\{
    \frac{1}{\cos (\wt - \theta)} \I_{\wt} (\lambda \cos(\wt  -\theta)^+)
    \,\Big| \,  \wt \in \ad (\overline{\D}) \right\}  & \,\, \mbox{if} \,\,
  \lambda \leq \eta_{\theta,h}\,,\\
+\infty & \,\, \mbox{if} \,\,  \lambda > \eta_{\theta,h}  \,.\end{array}
\right. $$
Moreover, if we define
$$ \delta_{\theta, h} \,=\, \delta \inf_{\wt \in \D} \frac{|\cos \wt| +
|\sin \wt|}{\cos (\wt - \theta)}\,. $$
where  $\delta = \inf \{\lambda \,|\, \PP (t(e) \leq \lambda) >0 \}$, the
good rate function $\K$ has the following properties: it is continuous on
$[0, \eta_{\theta,h}]$ except possibly at $\delta_{\theta,h}$ where it may
be only right continuous, it is infinite
 on $[0, \delta_{\theta, h} [ \cup ]\eta_{\theta,h}, +\infty[ $, finite on
 $]\delta_{\theta, h} ,\eta_{\theta,h} ]$, positive on $[\delta_{\theta, h}
 ,\eta_{\theta,h} [$ and equal to $0$ at $\eta_{\theta, h}$, and strictly
 decreasing when it is finite, in the sense that if $\K(\lambda) <\infty$,
 for all $\eps >0$, $\K(\lambda - \eps) > \K(\lambda)$.
\end{thm}

\begin{rem}
\label{chapitre6pgdcasdroit}
We will prove in Lemma \ref{chapitre6Kcasdroit}, section \ref{chapitre6propK} that when
$\theta \in \{0, \pi/2 \}$,
we have $\K(\lambda) = \I_{\theta} (\lambda ^+)$, and so Theorem \ref{thmPGDphi}
is consistent with the large deviation principle obtained in Theorem 3.17
by \cite{RossignolTheret08b} in the case of straight cylinders.
\end{rem}


\subsection{Comments on the hypotheses}
\label{remcond}

We want to discuss a little bit the different conditions on $F$ and $h$ we
use. The condition
$\mathbf{(H1)}$ is needed to obtain asymptotic results independent of the
height function $h$. The condition $\mathbf{(F2)}$ is needed to define $\nu_{\theta}$ in
the way we did it (it may be relaxed, see Remark 2.6 in
\cite{RossignolTheret08b}). The condition $\mathbf{(F1)}$ is equivalent
to the fact that $\nu_{\theta} \neq 0$.

The conditions $\mathbf{(F4)}$ or $\mathbf{(H3)}$ appear in Theorem
\ref{thmPGDphi} to deal with the upper bound of the large deviation principle
(see section \ref{seculd}, in particular Remark \ref{remcondsup}). They
correspond to the condition $\mathbf{(F5)}$ in Theorem
\ref{thmPGDtau}. Indeed, we need a stronger moment condition to deal with
the upper large deviations of $\tau$ because a minimal cutset corresponding to
the maximal flow $\tau(nA,h(n))$ is pinned along the boundary $\partial
(nA)$ of $nA$, thus it suffices that some edges in a neighbourhood of this
boundary have a huge capacity to increase the variable $\tau(nA,h(n))$ (for
more details, see section 2.3 in \cite{Theret:uppertau}). Since a minimal
cutset corresponding to the flow $\phi(nA,h(n))$ is not pinned, there are
no edges in the cylinders with such an influence on $\phi(nA,h(n))$.
However, the fact that a cutset for $\phi(nA,h(n))$ is not pinned implies
that it can be located anywhere in the cylinder $\cyl(nA,h(n))$, thus we
need to control the height of the cylinder by the condition $\mathbf{(H2)}$
to obtain interesting results concerning $\phi(nA,h(n))$, whereas this
condition does not appear in theorems concerning $\tau(nA,h(n))$. As
explained in Remark 3.22 in \cite{RossignolTheret08b}, the condition
$\mathbf{(H2)}$, combined with $\mathbf{(F1)}$, is relevant to observe a
maximal flow $\phi(nA,h(n))$ that is not null.

Finally, the conditions
$\mathbf{(FH1)}$ and $\mathbf{(FH2)}$ also appear because the minimal
cutset does not have fixed boundary conditions, thus it chooses its
orientation to solve an optimisation problem. The condition
$\mathbf{(FH1)}$ ensures that the direction chosen by an optimal cutset is
stable when $n$ goes to infinity; this condition, combined with
$\mathbf{(F2)}$ and $\mathbf{(H1)}$, is relevant to observe a
limit for $\phi(nA,h(n))/(nl(A))$, as proved in Theorem
\ref{chapitre6thm:lgn}. The condition $\mathbf{(FH2)}$ is of the same kind.

For simplicity of notation, we will denote by $\phi_n$ the maximal flow
$\phi(nA,h(n))$ for given $A$ and $h$ clearly given in the context. We shall often use two abbreviations: \emph{l\`agl\`ad} for "limite \`a gauche, limite \`a
droite", meaning that a function admits, on every point of its domain,
a limit (eventually infinite) from the left and a limit from the
right. We shall also use \emph{l.s.c} for "lower
semi-continuous".

\section{Lower large deviations}
\label{chapitre6subsec:deviations}

This section is devoted to the study of $\PP [\phi_n \leq \lambda n l(A)]$
for $\lambda \geq 0$. We will add conditions on $h$ and $F$ step
by step, to emphasize what condition is needed at each time.



\subsection{Technical lemma}
\label{chapitre6subsubsec:lemtech}

For any angle $\theta$, we define two vectors by their coordinates:
$$ \va \,=\, (\cos \theta, \sin \theta ) \quad \textrm{and} \quad \vc
\,=\, (\sin \theta, -\cos \theta)\,.  $$
We state here a property which comes from the weak triangle inequality for
$\nu$ (see section 4.4 in \cite{RossignolTheret08b}):

\begin{lem}
\label{chapitre6trigI}
Let $(abc)$ be a non degenerate triangle in $\RR^2$ and let $v_a$, $v_b$,
$v_c$ be the exterior normal unit vectors to the sides $[bc]$, $[ac]$,
$[ab]$. We denote by $(\cos \wt_i, \sin \wt_i)$ the coordinates of $v_i$,
and by $l(ij)$ the length of the side $[i,j]$ for $i$, $j$ in
$\{a,b,c\}$. If the angles $\widehat{cab}$ and $\widehat{abc}$ have values strictly
smaller than $\pi/2$, then for all $\lambda \geq 0$, for all $\alpha\in
[0,1]$, we have
$$ l(ab) \I_{\wt_c} \left(\frac{\lambda}{l(ab)} ^+ \right) \,\leq\, l(ac)
\I_{\wt_b} \left( \alpha \frac{\lambda}{l(ac)} ^+ \right) + l(bc) \I_{\wt_a}
\left((1-\alpha) \frac{\lambda}{l(bc)}^+ \right)\,.  $$
\end{lem}

\begin{dem}
This proof follows the one of proposition 11.6 in \cite{Cerf:StFlour}. We
consider the cylinder
$$ \cyl_c(N) \,=\, \cyl (N[ab], N) $$
of dimensions $N l(ab)\times 2N$ oriented towards the direction $\wt_c$,
and we define $\tau_c(N) = \tau(\cyl_c(N))$ (implicitly, for the direction
defined by $\wt_c$). Exactly as
in section 4.1 of \cite{RossignolTheret09}, we choose two functions $\zeta, h': \NN
\rightarrow \RR^+$ such that
$$ \lim_{n\rightarrow \infty} h'(n) \,=\, \lim_{n\rightarrow \infty}
\zeta(n) \,=\, + \infty \,, $$
and
$$\lim_{n\rightarrow \infty} \frac{h'(n)}{\zeta(n)} \,=\, 0\,. $$
We construct smaller cylinders oriented towards the directions $\wt_b$ and
$\wt_a$ inside $\cyl_c(N)$ (see figure \ref{chapitre6emboitement3}).
\begin{figure}[!ht]
\centering
\input{emboitement3.pstex_t}
\caption{The cylinders $\cyl_c(N)$, $\cyl^i_b(n)$ and $\cyl^j_a(n)$.}
\label{chapitre6emboitement3}
\end{figure}
We define
$$\cyl_b(n) \,=\, \cyl \left( [0, 0+nl(ab) \vg] ,  h'(n) \right) \,,$$
$$ \cyl_a(n) \,=\, \cyl \left([0,0+nl(bc) \vh] ,  h'(n) \right) \,,$$
respectively oriented towards the direction $\wt_b$ and $\wt_a $. We define
the vectors
$$ \vec{u}_i \,=\, \left( \zeta(n) + (i-1) nl(ac) \right) \vg \quad
\textrm{and} \quad  \vec{w}_j \,=\, \left( \zeta(n) + (j-1)nl(bc) \right)
\vh  $$
and the points
$$ U_i \,=\, Na + \vec{u}_i \quad \textrm{and} \quad W_j \,=\, Nc + \vec{w}_j $$
for 
$$i\,\in\, \left\{ 1,..., \M_b=\left\lfloor \frac{Nl(ac) -2\zeta(n)}{nl(ac)}
  \right\rfloor  \right\}\quad \textrm{and} \quad j\,\in\, \left\{ 1,...,
  \M_a=\left\lfloor \frac{Nl(bc) -2\zeta(n)}{nl(bc)} \right\rfloor  \right\} \,,$$
where $\M_a=\M_a(n,N,b,c)$ and $\M_b=\M_b(n,N,a,c)$. For $i=1,...,\M_b$
(resp. $j=1,...,\M_a$), let $\widetilde{\cyl}_b^i(n)$ (resp. $\widetilde{\cyl}_a^j(n)$) be the
image of $\cyl_b(n)$ (resp. $\cyl_a(n)$) the the translation of vector
$\overrightarrow{0U_i}$ (resp. $\overrightarrow{0W_j}$). We can translate
again each $\widetilde{\cyl}_b^i(n)$ (resp. $\widetilde{\cyl}_a^j(n)$) by a
vector of norm strictly smaller than one to obtain a integer translate
$\cyl_b^i(n)$ (resp. $\cyl_a^j(n)$) of $\cyl_b(n)$ (resp. $\cyl_a(n)$),
i.e., a translate by a vector whose coordinates are in $\mathbb{Z}^2$. For
$i=1,...,\M_b$ (resp. $j=1,...,\M_a$), we define $\tau^i_b(n) =
\tau (\cyl^i_b(n))$ (resp. $\tau_a^j(n) = \tau (\cyl_a^j(n))$) for the
direction  defined by $\wt_b$ (resp. $\wt_a$). The dimensions of $\cyl_b(n)$
(resp. $\cyl_a(n)$) are $(n l(ac)) \times
2h'(n)$ (resp. $(n l(bc)) \times
2h'(n)$), and for $N$ and $n$ large enough $\cyl^i_b(n)$ and $\cyl^j_a(n)$ are included
in $\cyl_c(N)$ for all $i$ and $j$  (we only consider such large $n$ and
$N$), because $\widehat{cab}$ and
$\widehat{abc}$ are strictly smaller than $\pi /2$ and $h'(n)/\zeta(n)
\rightarrow 0$. The variables $(\tau_b^i(n), \tau_a^j(n))$ are identically
distributed. To glue together cutsets in the cylinders $\cyl^i_b(n)$ and
$\cyl^j_a(n)$ for all $i$ and $j$ to obtain a cutset in $\cyl_c(N)$
we have to add some edges. We finally define, for a constant $\zeta \geq 4$,
$$ \E_3 (n,N,a,b,c) =\mathcal{V} \left( 
\begin{array}{l}
[Na,Na +\vec{u}_1] \cup [Na +\vec{u}_{\M_b}, Nc]\\
 \cup
  [Nc, N c + \vec{w}_1] \cup [ Nc +\vec{w}_{\M_a} ,Nb] \end{array} ,\zeta \right) ,$$
and we denote by $E_3(n,N,a,b,c)$ the set of the edges included in
$\E_{3}(n,N,a,b,c)$. There exists a constant $C_7$ such that
$$ \card (E_3 (n,N, a,b,c)) \,\leq\, C_7 \left(\zeta(n) + n + \frac{N}{n} \right) \,. $$
The union of $E_3(n,N,a,b,c)$ with cutsets in the cylinders $\cyl_b^i(n)$
and $\cyl_a^j(n)$ for all $i$ and $j$ separates the upper half part from
the lower half part of the boundary of $\cyl_c(N)$ (see figure
\ref{chapitre6emboitement3}), so we have
\begin{equation}
\label{chapitre6eqbase}
\tau_c(N) \,\leq\, \sum_{i=1}^{\M_b}\tau^i_b(n) + \sum_{j=1}^{\M_a}\tau^j_a(n) + V(E_3(n,N,a,b,c))\,.
\end{equation}
Then for all $\lambda \geq 0$, for all positive $\eta$, for all large $N$,
for all $\alpha \in [0,1]$, by the FKG inequality we have
\begin{align*}
\PP \Bigg[ \frac{\tau_c(N)}{Nl(ab)} & \leq  \lambda + 3\eta -
\frac{1}{\sqrt{Nl(ab)}} \Bigg]\\ & \,\geq\, \PP \Bigg[ \frac{\tau_c(N)}{Nl(ab)} \leq \lambda + 2\eta
\Bigg]\\
& \,\geq\, \PP \left[\sum_{i=1}^{\M_b} \tau^i_b(n) \leq \alpha (\lambda +\eta) N l(ab) \right] \times \PP \left[ \sum_{j=1}^{\M_a}\tau^j_a(n) \leq
  (1-\alpha) (\lambda +\eta) N l(ab) \right] \\ & \qquad \times \PP \left[ V( E_3(n,N,a,b,c)) \leq \eta N l(ab)  \right]\\
& \,\geq\, \prod_{i=1}^{\M_b} \PP \left[ \tau^i_b(n) \leq \alpha (\lambda +\eta) n l(ab) \right] \times \prod_{j=1}^{\M_a} \PP \left[\tau^j_a(n) \leq
  (1-\alpha) (\lambda +\eta) n l(ab) \right]\\ & \qquad \times \PP \left[ V( E_3(n,N,a,b,c)) \leq \eta N l(ab)  \right]\\
&  \,\geq\, \PP \left[ \tau_a(n) \leq (1-\alpha) (\lambda +\eta) n l(ab) -
  \frac{1}{\sqrt{nl(bc)}} \right]\\ &\hspace{-2cm}\times \PP \left[ \tau_b(n) \leq \alpha (\lambda +\eta) n l(ab) - \frac{1}{\sqrt{nl(ac)}} \right]  \times \PP \left[ t(e) \leq \frac{\eta N l(ab)}{C_7( \zeta(n)+n+N/n)}
\right]^{C_7 (\zeta(n)+n+N/n )} \,.
\end{align*}
We take the logarithm of the previous inequality, divide it by $-N$, send
$N$ to infinity and then $n$ to infinity. We obtain that
\begin{equation}
\label{chapitre6base}
l(ab) \I_{\wt_c} \left( \lambda + 3\eta \right) \,\leq\, l(ac)
\I_{\wt_b} \left( \alpha (\lambda + \eta)\frac{l(ab)}{l(ac)} \right) + l(bc) \I_{\wt_a}
\left( (1-\alpha) (\lambda + \eta) \frac{l(ab)}{l(bc)} \right) \,.
\end{equation}
Sending $\eta$ to zero, we obtain the desired inequality.
\end{dem}

We state next a property of continuity:
\begin{lem}
\label{chapitre6lem1}
For all $\lambda \geq 0$, we define $g_{\lambda}:[\theta-\pi/2,
\theta+\pi/2] \rightarrow \RR^+ \cup \{+\infty\}$ by 
$$ \forall \wt \in ]\theta - \pi/2, \theta + \pi/2[, \qquad g_{\lambda}(\wt) \,=\,
\frac{1}{\cos (\wt - \theta)} \I_{\wt} (\lambda \cos(\wt - \theta)^+) $$
and
$$ g_{\lambda}(\theta-\pi/2) \,=\, g_{\lambda}
(\theta+\pi/2)\,=\,\left\{  \begin{array}{ll} +\infty & \textrm{if } \I_{\wt}(0^+) >0
    \,, \\ 0 & \textrm{if } \I_{\wt}(0^+)=0  \,. \end{array} \right. $$
Then $g_{\lambda}$ is lower semi-continuous, and $g_{\lambda}$ is continuous on 
$$ H^{>}_\lambda \,=\, \left\{ \wt \,|\, \lambda > \delta
  \frac{|\cos \wt| + |\sin \wt|}{\cos (\wt - \theta)} \right\}\,. $$
\end{lem}
\begin{rem}
\label{chapitre6reminfini}
If $\mathbf{(F1)}$ holds, then for all $\wt$ we have $\nu_{\wt} >0$, that
implies $\I_{\wt} (0^+) >0$. The definition of $ g_{\lambda}(\theta-\pi/2)$
and $ g_{\lambda}(\theta+\pi/2)$ is consistent with the expression given
for any different $\wt$. We shall always use Lemma \ref{chapitre6lem1}
under assumption $\mathbf{(F1)}$.
\end{rem}

\begin{dem}
The proof is based on the same ideas as the one of lemma \ref{chapitre6trigI}, so
we will use part of it. We consider two angles $\wt_1$, $\wt_2$ such that
$\wt_1 - \wt_2 = \hat{\varepsilon}$ (positive or negative) and
$|\hat{\varepsilon} | = \varepsilon$ is small. Let $(abc)$ be the right
triangle such that, using the same notations as in the previous proof,
$l(ab)=1$, $\wt_c = \wt_1 + \pi$, $\wt_b=\wt_2$ and $\wt_a = \wt_2 -
\pi/2$, and so $\widehat{bac} = \eps$, $\widehat{acb} = \pi/2$ and $\widehat{abc}<\pi/2$.
Obviously we are confronted with a particular case of triangle $(abc)$
studied in lemma \ref{chapitre6trigI}. We do exactly the same construction as in the
previous proof, and we start again from equation (\ref{chapitre6base}). Here we have
constructed $(abc)$ such that $l(ab)=1$, $l(ac) = \cos \eps$ and
$l(bc)=\sin \eps$, and by invariance of the graph by a rotation of angle
$\pi/2$, we know that the functions $\I_{\wt_2}$ and $\I_{\wt_2 - \pi/2}$
(respectively $\I_{\wt_1}$ and $\I_{\wt_1 + \pi}$)
are equal. We can rewrite equation (\ref{chapitre6base}) the following way:
\begin{equation}
\label{chapitre6base2}
\I_{\wt_1} (\lambda + 3\eta) \,\leq\, (\cos \eps) \I_{\wt_2} \left( \alpha
  \frac{\lambda + \eta}{\cos \eps} \right) + (\sin \eps) \I_{\wt_2} \left(
  (1-\alpha) \frac{\lambda + \eta}{\sin \eps} \right)\,.
\end{equation}
We want to make appear the factor $\cos (\wt_1 - \theta)$, so for all
$\lambda \geq 0$ and for all small $\eta$ we deduce from (\ref{chapitre6base2}) that
for all $\eps$ small enough,
\begin{align*}
\I_{\wt_1} &(\lambda \cos(\wt_1 - \theta) + 3\eta)\\
& \,\leq\, (\cos \eps) \I_{\wt_2}
\left( \alpha \frac{\lambda \cos(\wt_1-\theta) + \eta}{\cos \eps} \right) +
(\sin \eps) \I_{\wt_2} \left( (1-\alpha) \frac{\lambda \cos(\wt_1 - \theta)
  + \eta}{\sin \eps} \right)\\
& \,\leq\, (\cos \eps ) \I_{\wt_2} \left( \alpha (\lambda \cos(\wt_2
  -\theta) + \eta/2) \right) + (\sin \eps) \I_{\wt_2} \left(
  (1-\alpha)\frac{\lambda \cos(\wt_1 - \theta) + \eta}{\sin \eps}  \right)\,.
\end{align*}
If $\lambda >0$ we choose $\alpha \in] \max(2/3 , 1-\eta/(12 \lambda)) ,1[$
(remember that $\lambda$ is fixed and we can choose $\eta$ small in
comparison with $\lambda$), then $\alpha (\lambda \cos (\wt_2 -\theta) +
\eta/2) \geq \lambda \cos(\wt_2 - \theta) + \eta /4$. This equation is
satisfied for all $1>\alpha \geq 1/2$ if $\lambda =0$. We stress here the
fact that how large must be $\alpha$ depends on $\lambda$ and $\eta$, but
not on $\eps$. With a such fixed big $\alpha$, we obtain that
\begin{align*}
\I_{\wt_1} &(\lambda \cos(\wt_1 - \theta) + 3\eta) \\
&  \,\leq\, (\cos \eps ) \I_{\wt_2} \left( \lambda \cos(\wt_2
  -\theta) + \eta/4 \right) + (\sin \eps) \I_{\wt_2} \left(
  (1-\alpha)\frac{\lambda \cos(\wt_1 - \theta) + \eta}{\sin \eps}  \right)\,.
\end{align*}
We send $\wt_2$ to $\wt_1$, i.e. $\eps$ to zero by fixing $\wt_1$. Since $(1-\alpha)(\lambda \cos(\wt_1 -
\theta) + \eta)$ is fixed and positive, we know that for small $\eps$ we obtain
$$ (1-\alpha)\frac{\lambda \cos(\wt_1 - \theta) + \eta}{\sin \eps} \,>\,
\nu_{\max} \,=\, \max_{\theta \in [0,\pi]} \nu_{\theta}\,, $$
and so for all $\wt$ we have
$$ \I_{\wt} \left(
  (1-\alpha)\frac{\lambda \cos(\wt_1 - \theta) + \eta}{\sin \eps}
\right)\,=\,0\,. $$
We send finally $\eta$ to zero and obtain
$$\I_{\wt_1} (\lambda \cos (\wt_1 - \theta) ^+) \,\leq\, \liminf_{\eta
  \rightarrow 0} \liminf_{\hat{\varepsilon} \rightarrow 0} \I_{\wt_1 +
  \hat{\varepsilon}} \left( \lambda \cos (\wt_1 + \hat{\varepsilon} - \theta) +
  \eta/4 \right) \,. $$
We know that the limit $\lim_{\eta
  \rightarrow 0} \I_{\wt_1 + \hat{\varepsilon}} (\lambda \cos (\wt_1 +
\hat{\varepsilon} - \theta) + \eta/4)$ is an increasing limit for all fixed
$\hat{\varepsilon}$, so we get:
\begin{equation}
\label{chapitre6ccl1}
\I_{\wt_1} (\lambda \cos (\wt_1 - \theta ) ^+)\,\leq\, \liminf_{\hat{\varepsilon}
  \rightarrow 0} \I_{\wt_1 + \hat{\varepsilon} } (\lambda \cos(\wt_1 +
\hat{\varepsilon} - \theta) ^+) \,.
\end{equation}

We will now fix $\wt_2$ and send $\wt_1$ to $\wt_2$. Starting again from
(\ref{chapitre6base}), for all $\beta >0$, for all $\lambda >0$, for all $\wt_2 \in
]\theta - \pi/2, \theta + \pi/2[$, for all $\eta$ small enough and $\eps$
small (in particular such that
$\wt_1 \in]\theta - \pi/2, \theta + \pi/2[ $ too), we obtain
\begin{align*}
\I_{\wt_1}  &(\lambda \cos(\wt_1 - \theta) + \beta)\\
& \,\leq\, \I_{\wt_1} (\lambda \cos(\wt_1 - \theta))\\
& \,\leq\, (\cos \eps) \I_{\wt_2} \left(\alpha \frac{\lambda \cos (\wt_1
    -\theta) -2\eta}{\cos \eps}  \right) + (\sin \eps) \I_{\wt_2} \left(
  (1-\alpha) \frac{\lambda \cos (\wt_1 -\theta) -2 \eta}{\sin \eps}
\right)\\
& \,\leq\, (\cos \eps) \I_{\wt_2} \left(\alpha( \lambda \cos (\wt_2 -
  \theta)- 3\eta) \right) + (\sin \eps) \I_{\wt_2} \left(
  (1-\alpha) \frac{\lambda \cos (\wt_1 -\theta) - 2\eta}{\sin \eps}
\right)\,.
\end{align*}
Exactly as previously, for $\alpha <1$ but sufficiently close to $1$ (how
close depending on $\lambda$ and $\eta$ but not on $\eps$), we have
\begin{align*}
\I_{\wt_1}  &(\lambda \cos(\wt_1 - \theta) + \beta)\\
& \,\leq\, (\cos \eps) \I_{\wt_2} \left(\lambda \cos (\wt_2 -
  \eps)- 4\eta \right) + (\sin \eps) \I_{\wt_2} \left(
  (1-\alpha) \frac{\lambda \cos (\wt_1 -\theta) - 2\eta}{\sin \eps}
\right)\,.
\end{align*}
We send first $\beta$ to zero, then $\wt_1$ to $\wt_2$ (thus $\eps $ to
zero), and finally $\eta$ to zero to obtain as for (\ref{chapitre6ccl1}) that
\begin{equation}
\label{chapitre6ccl2}
\I_{\wt_2} (\lambda \cos(\wt_2 - \theta) ^-) \,\geq\, \limsup_{\hat{\varepsilon}
  \rightarrow 0} \I_{\wt_2 + \hat{\varepsilon}} (\lambda \cos (\wt_2 +
\hat{\varepsilon} - \theta)^+)\,.
\end{equation}
This inequality remains valid for $\lambda = 0$ or $\cos(\wt_2 -
\theta)=0$, since for convenience we decided that $\I_{\wt_2}(0^-) = +\infty$.
From (\ref{chapitre6ccl1}) and (\ref{chapitre6ccl2}), we conclude that for all $\lambda \geq 0$:
\begin{align*}
 \frac{1}{\cos(\wt -\theta)}\I_{\wt} (\lambda \cos (\wt - \theta )^+)
 \,\leq\, \ & \liminf_{\hat{\varepsilon}
  \rightarrow 0} g_{\lambda}(\wt + \hat{\varepsilon}) \,\leq\,\\
& \,\leq\,  \limsup_{\hat{\varepsilon}
  \rightarrow 0} g_{\lambda}(\wt + \hat{\varepsilon}) \,\leq\,\frac{1}{\cos(\wt -\theta)} \I_{\wt}
(\lambda \cos (\wt - \theta) ^-) \,.  
\end{align*}
Lemma \ref{chapitre6lem1} follows, since we know that:
$$\forall \wt \in H^{>}_{\lambda} \qquad \I_{\wt}(\lambda \cos(\wt -
\theta)^+) \,=\, \I_{\wt} (\lambda \cos (\wt - \theta)^- ) \,.$$
\end{dem}


\subsection{Lower bound}
\label{chapitre6subsubsec:lowerboundpgd}

From now on, we suppose that the height function $h$ satisfies $\mathbf{(H1)}$.
We will use equation (19) of \cite{RossignolTheret09}, and thus the
construction that leads to it, to prove that
\begin{equation}
\label{chapitre6eq1'}
\liminf_{n\rightarrow \infty } \frac{1}{nl(A)} \log \PP [\phi_n \leq
\lambda n l(A)] \,\geq\, -  \inf_{\wt \in \ad(\underline{\D})}
\frac{1}{\cos (\wt - \theta)} \I_{\wt} \left(\lambda \cos (\wt - \theta) ^-
\right)\,,
\end{equation}
and
\begin{equation}
\label{chapitre6eq1bis'}
\limsup_{n\rightarrow \infty } \frac{1}{nl(A)} \log \PP [\phi_n \leq
\lambda n l(A)] \,\geq\, -  \inf_{\wt \in \ad(\overline{\D})}
\frac{1}{\cos (\wt - \theta)} \I_{\wt} \left(\lambda \cos (\wt - \theta) ^-
\right)\,.
\end{equation}
We recall this construction here. We
consider a line segment $A$, of
orthogonal unit vector $\va = (\cos\theta, \sin\theta)$ for $\theta
\in [0, \pi[$, and a function $h: \NN \rightarrow \RR^+ $ satisfying $\lim_{n\rightarrow
  \infty}h(n) = +\infty$. We use the notation $\D_n=\D(nA,h(n))$. For all $\wt \in
\D_n $,  we define
$$  k_n \,=\, \frac{1}{2} + \frac{nl(A) \tan (\wt - \theta)}{4 h(n)} \,,$$
and thus $\kappa_n=(k_n,\widetilde{\theta})\in D_n$. We want to
compare $\phi_n^{\kappa_n}$ with the maximal flow $\tau$ in a cylinder
inside $\cyl(nA, h(n))$ and oriented towards the direction
$\widetilde{\theta}$. In fact, we must use the subadditivity of $\tau$ and
compare $\phi_n^{\kappa_n}$ with a sum of such variables $\tau$.

\begin{figure}[!ht]
\centering
\input{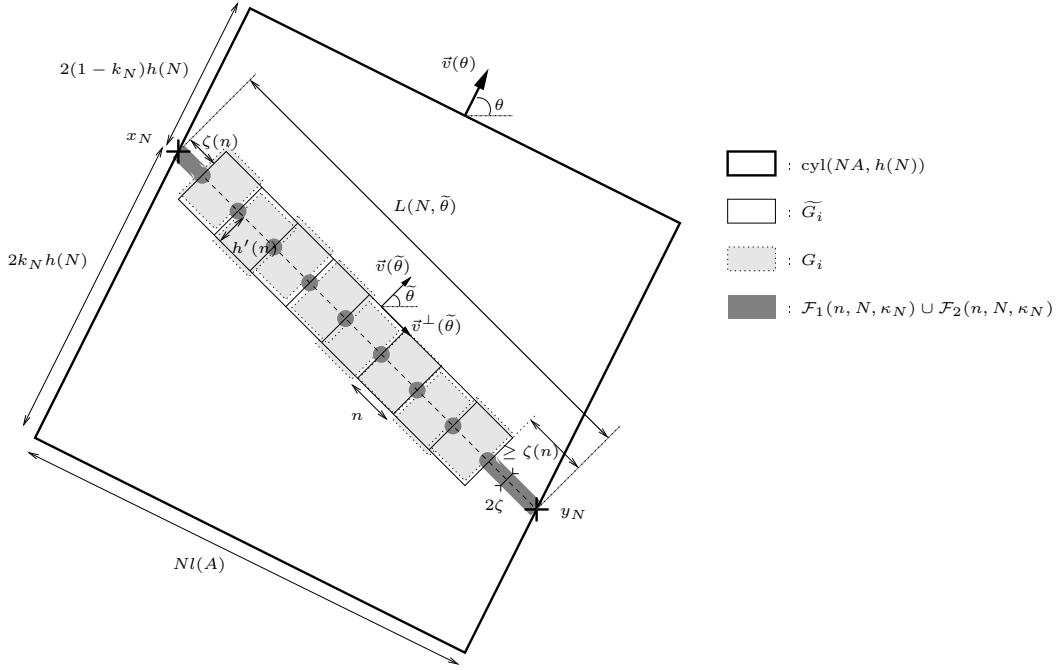}
\caption{The cylinders $\cyl(NA,h(N))$ and $G_i$, for $i=1,...,\M$.}
\label{chapitre6emboitement1}
\end{figure}

We consider $n$ and $N$ in $\NN$, with $N$ a lot bigger than $n$. The following definitions can seem a little bit
complicated, but Figure \ref{chapitre6emboitement1} is more explicit.
We choose two functions $h', \zeta :\NN \rightarrow \RR^+$ such that
$$\lim_{n\rightarrow \infty} h'(n) \,=\, \lim_{n\rightarrow \infty}
\zeta(n) \,=\, + \infty \,,$$
and 
\begin{equation}
\label{chapitre6cond1}
 \lim_{n\rightarrow \infty}\frac{h'(n)}{\zeta (n)} \,=\, 0 \,. 
\end{equation}
We consider a fixed $\wt \in \D_N$. We recall that
$$\vb \,=\, (\cos \widetilde{\theta}, \sin \widetilde{\theta}) \quad
\textrm{and} \quad \vd \,=\, (\sin
\widetilde{\theta}, -\cos \widetilde{\theta}) \,.$$
In $\cyl(NA,h(N))$, we denote by $x_N$ and $y_N$ the two points
corresponding to the boundary conditions $\kappa_N$, such that
$\overrightarrow{x_Ny_N} \cdot \vd >0 $. Notice that according to our
choice of $k_N$, the segments $[x_N,y_N]$ and $NA$ cut each other in their middle. If we denote by $L(N,\wt)$ the
distance between $x_N$ and $y_N$, we have:
$$ L(N,\wt)\,=\, \frac{Nl(A)}{\cos(\widetilde{\theta}-\theta)} \,. $$
We define
$$ \cyl'(n) \,=\, \cyl ([0 ,n \vd ],h'(n))  \,.$$
We will translate $\cyl '(n)$ numerous times inside $\cyl(NA,h(N))$. We
define
$$ t_i \,=\, x_N + \left(\zeta(n) + (i-1) n \right) \vd \,,$$
for $i=1,..., \M$, where
$$ \M \,=\, \M(n,N) \,=\, \left\lfloor \frac{L(N,\wt)  -2\zeta(n)}{n} \right\rfloor \,.$$
Of course we consider only $N$ large enough to have $\M \geq 2$. For
$i=1,...,\M$, we denote by $\widetilde{G_i}$ the image of $\cyl'(n)$ by the
translation of vector $\overrightarrow{0 t_i}$. For $n$ (and thus $N$)
sufficiently large, thanks to condition (\ref{chapitre6cond1}), we know that
$\widetilde{G_i} \subset \cyl(NA,h(N))$ for all $i$. We can translate $\widetilde{G_i}$ again by a vector of
norm strictly smaller than $1$ to obtain an integer translate of $\cyl'(n)$
(i.e., a translate by a vector whose coordinates are in $\ZZ^2$) that we
will call $G_i$. Now we want to glue together cutsets of boundary condition
$(1/2, \wt)$ in the cylinders $G_i$. We define:
$$ \F_1(n,N,\kappa_N) \,=\, \left(\bigcup_{i=1}^{\M} \mathcal{V} (t_i,
  \zeta_0) \right) \, \bigcap \,\cyl(NA,h(N))  \,,$$
where $\zeta_0$ is a fixed constant larger than $4$, and:
$$ \F_2(n,N,\kappa_N) \,=\, \mathcal{V} \left( [x_N,x_N + \zeta(n) \vd ]
  \cup [z_\M , y_N ] ,\zeta_0 \right)\, \bigcap \,\cyl(NA,h(N)) \,.$$
Let $F_1(n,N,\kappa_N)$ (respectively $F_2(n,N,\kappa_N)$) be the set of
the edges included in $\F_1(n,N,\kappa_N)$ (respectively
$\F_2(n,N,\kappa_N)$). If for every $i=1,...,\M$, $\G_i$ is a cutset of
boundary condition $(1/2,\wt)$ in $G_i$, then
$$ \bigcup_{i=1}^{\M} \G_i \cup F_1(n,N,\kappa_N) \cup F_2(n,N,\kappa_N) $$
contains a cutset of boundary conditions $\kappa_N$ in $\cyl(NA,h(N))$. We
obtain:
\begin{equation}
\label{chapitre6lien_phiF_tau_1} 
\phi_N^{\kappa_N} \,\leq\, \sum_{i=1}^{\M} \tau(G_i,\vb ) +
V(F_1(n,N,\kappa_N) \cup F_2(n,N,\kappa_N) ) \,,
\end{equation}
and so,
\begin{equation}
\label{chapitre6lien_phi_tau} 
\forall \widetilde{\theta}\in \D_N \qquad \phi_N \,\leq\,
\phi_{N}^{\kappa_N} \,\leq\, \sum_{i=1}^{\M}\tau(G_i, \vb ) + V( F_1(n,N,\kappa_N) \cup F_2(n,N,\kappa_N) ) \,.
\end{equation}
This equation (\ref{chapitre6lien_phi_tau}) is equation (19) in
\cite{RossignolTheret09}. Moreover, there exists a constant $C_5$ such that:
$$ \card ( F_1(n,N,\kappa_N) ) \,\leq\, C_5
  \M  \quad \textrm{and} \quad \card (F_2(n,N,\kappa_N))\,\leq\,  C_5 \left( \zeta(n) + n \right) \,.$$
Then for all $\widetilde{\theta}\in \D_n$, for all $\lambda >0$,
for all positive small $\varepsilon$, by the FKG inequality,
\begin{align*}
\PP [  \phi_N \leq \lambda l(A) N] & \,\geq\, \PP \left[
\bigcap_{i=1}^{\M} \left\{\tau(G_i, \vb ) \leq (\lambda - \varepsilon)
  \frac{Nl(A)}{\M} \}\cap
\{ V(E(n,N,\kappa_N)) \leq \varepsilon l(A) N  \right\}\right]\\
& \,\geq\, \PP \left[
 \tau(\cyl'(n),\vb )  \leq (\lambda - \varepsilon) \cos(\wt - \theta) n
\right] ^{\M}\\ & \qquad \qquad \times \PP
\left[ \forall e\in E(n,N,\kappa_N) \,,\, t(e) \leq \frac{ \varepsilon
    l(A) N }{C_5 (\M + \zeta(n) + n)} \right]\\
& \,\geq\, \PP \left[
  \frac{\tau(\cyl'(n),\vb )}{n}  \leq (\lambda -\varepsilon) \cos
  (\widetilde{\theta}  -\theta) - \frac{1}{\sqrt{n}} \right]^{\M}\\ &\qquad \qquad \times \PP \left[ t(e) \leq \frac{ \varepsilon
    l(A) N }{C_5 (\M+\zeta(n)+n)} \right]^{C_5 (\M + \zeta(n) +n)}\,.
\end{align*}
We take the logarithm of this inequality, divide it by $Nl(A)$, send $N$
to infinity and then $n$ to infinity. Thanks to Theorem \ref{propI},
for all $\widetilde{\theta} \in \underline{\D}$ and $\lambda >\varepsilon >0$, we obtain
$$ \liminf_{N\rightarrow \infty }
\frac{1}{Nl(A)} \log \PP \left[\frac{\phi_N}{Nl(A)} \leq \lambda \right] \,\geq\, \frac{-1}{\cos
(\widetilde{\theta} - \theta)} \I_{\widetilde{\theta}}\left(
(\lambda-\varepsilon) \cos( \widetilde{\theta} - \theta)  \right) \,. $$
Sending $\varepsilon$ to zero (remember that $\I_{\wt}$ is l\`agl\`ad) and
taking the infimum in $\widetilde{\theta}$,
\begin{equation}
\label{chapitre6eq1}
\liminf_{N\rightarrow \infty }
\frac{1}{Nl(A)} \log \PP [\phi_N \leq \lambda Nl(A)] \,\geq\, -
\inf_{\widetilde{\theta} \in \underline{\D}}  \frac{1}{\cos
  (\widetilde{\theta} - \theta)} \I_{\widetilde{\theta}}\left( \lambda
  \cos( \widetilde{\theta} - \theta) ^- \right)\,.
\end{equation}
Similarly, if $\wt \in \overline{\D}$, let $\psi:\NN\rightarrow \NN$ be
strictly increasing such that for all $N$, $\wt \in \D_{\psi(N)}$. Then we
obtain by the same arguments that
\begin{align}
\label{chapitre6eq1bis}
 \limsup_{N\rightarrow \infty }
\frac{1}{Nl(A)} \log \PP [\phi_N \leq \lambda Nl(A)]& \,\geq\,
 \liminf_{N\rightarrow \infty }
\frac{1}{\psi(N)l(A)} \log \PP [\phi_{\psi(N)} \leq \lambda
\psi(N)l(A)]\nonumber \\
& \,\geq\, - \inf_{\widetilde{\theta} \in \overline{\D}}  \frac{1}{\cos
  (\widetilde{\theta} - \theta)} \I_{\widetilde{\theta}}\left( \lambda
  \cos( \widetilde{\theta} - \theta) ^- \right)\;.
\end{align}
These inequalities remain valid for $\lambda = 0$, since $\I_{\wt}(0^-)
=+\infty$, so equations (\ref{chapitre6eq1}) and (\ref{chapitre6eq1bis}) are satisfied for all
$\lambda \geq 0$.

We will transform a little bit inequalities (\ref{chapitre6eq1}) and (\ref{chapitre6eq1bis})
to make it more useful for us in the proof of the large deviation principle
below. Actually, let us prove that:
\begin{equation}
\label{chapitre6bornes}
\inf_{\wt \in \D} \frac{1}{\cos (\wt - \theta)} \I_{\wt}
\left(\lambda \cos (\wt - \theta) ^- \right) \,=\, \inf_{\wt \in \ad
  (\D)} \frac{1}{\cos (\wt - \theta)} \I_{\wt} \left(\lambda
  \cos (\wt - \theta) ^- \right)\,,
\end{equation}
where $\D$ is an interval of $[\theta - \pi/2, \theta + \pi/2]$ which
  is centered at $\theta$ and symmetric with respect to $\theta$
(representing $\overline{\D}$ or $\underline{\D}$ here).
As we did previously, we define
$$ H_\lambda^* \,=\, \left\{\wt \,|\, \lambda * \delta \frac{|\cos \wt | + |\sin
  \wt|}{\cos (\wt - \theta)}  \right\} \,,$$
where $*$ represents $ <,\,>,\,\leq, \,\geq$ or $=$, and for simplicity of
notations we define also:
$$ \widetilde{g}_{\lambda}(\wt) \,=\, \frac{1}{\cos (\wt - \theta)} \I_{\wt}
\left(\lambda \cos (\wt - \theta) ^- \right)\,. $$
The function $\widetilde{g}_{\lambda}$ is
infinite on $H^{\leq}_\lambda$, and finite, continuous and equal to
$g_\lambda$ on $H^>_\lambda$. If $\D$ is included in
$H^{\leq}_\lambda$, then $\ad(\D)$ too because
$H^{\leq}_\lambda$ is closed, and then:
$$\inf_{\D} \widetilde{g}_\lambda \,=\, +\infty
\,=\,\inf_{\ad(\D)} \widetilde{g}_\lambda  \,.$$ 
Otherwise, $\D\cap H^>_\lambda$ is non empty, so
$\inf_{\D}  \widetilde{g}_\lambda$ is finite. If
$\ad(\D) \neq \D$ (otherwise the result is
obvious),  then $\D$ is open since it is symmetric with respect to
  $\theta$, and we denote by $\wt_1 $ and $\wt_2$ the two points of
$\ad(\D) \smallsetminus \D$. Either $
\widetilde{g}_\lambda$ is continuous at $\wt_1$ (respectively $\wt_2$), or $
\widetilde{g}_\lambda(\wt_1)$ (respectively $
\widetilde{g}_\lambda(\wt_2)$) is infinite, so
$$ \inf_{\D} \widetilde{g}_\lambda \,=\, \inf_{\ad(\D)}
\widetilde{g}_\lambda  \,, $$
and equation (\ref{chapitre6bornes}) is proved. Inequalities (\ref{chapitre6eq1}) and
(\ref{chapitre6eq1bis}) are equivalent to (\ref{chapitre6eq1'}) and (\ref{chapitre6eq1bis'}).


\subsection{Upper bound}
\label{chapitre6subsubsec:upperboundpgd}

\begin{figure}[!ht]
\centering
\input{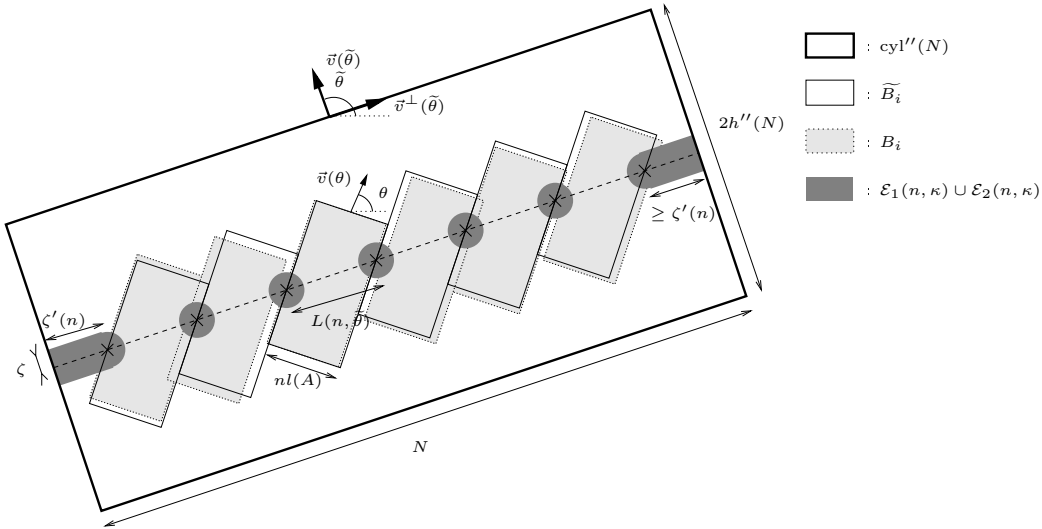}
\caption{The cylinders $\cyl''(N)$ and $B_i$, for $i=1,...,\N$.}
\label{chapitre6emboitement2}
\end{figure}

We suppose from now on that $h$ satisfies $\mathbf{(H1)}$ and $\mathbf{(H2)}$.
We will use equation (24) in \cite{RossignolTheret09}, and thus the
construction that leads to it, to prove that
\begin{equation}
\label{chapitre6eq2}
 \limsup_{n\rightarrow \infty }
\frac{1}{nl(A)} \log \PP [\phi_n \leq \lambda n l(A)] \,\leq\, -
\inf_{\widetilde{\theta} \in \ad (\overline{\D})}  \frac{1}{\cos
  (\widetilde{\theta} - \theta)} \I_{\widetilde{\theta}}\left( \lambda
  \cos( \widetilde{\theta} - \theta) ^+ \right)\,,
\end{equation}
and
\begin{equation}
\label{chapitre6eq2bis}
 \liminf_{n\rightarrow \infty }
\frac{1}{nl(A)} \log \PP [\phi_n \leq \lambda n l(A)] \,\leq\, -
\inf_{\widetilde{\theta} \in \ad (\underline{\D})}  \frac{1}{\cos
  (\widetilde{\theta} - \theta)} \I_{\widetilde{\theta}}\left( \lambda
  \cos( \widetilde{\theta} - \theta) ^+ \right)\,.
\end{equation}
We recall this construction now. We do the symmetric construction of the one done in section
\ref{chapitre6subsubsec:lowerboundpgd}. We consider $n$ and $N$ in $\NN$ and take $N$ a lot bigger than $n$. We
choose functions $\zeta', h'': \NN \rightarrow \RR^+$ such that
$$ \lim_{n\rightarrow \infty} \zeta'(n) \,=\, \lim_{n\rightarrow \infty}
h''(n) \,=\, +\infty \,, $$
and
\begin{equation}
\label{chapitre6cond2}
 \lim_{n\rightarrow \infty}\frac{ h(n)}{\zeta' (n)} \,=\,0 \,. 
\end{equation}
We consider $\kappa=(k,\widetilde{\theta})\in D_n$. Keeping the same
notations as in section \ref{chapitre6subsubsec:lowerboundpgd}, we define
$$ \cyl''(N) \,=\, \cyl \left([0,N \vd ] ,
  h''(N)\right) \,. $$
We will translate
$\cyl (nA, h(n))$ numerous times in $\cyl ''(N)$. The figure
\ref{chapitre6emboitement2} is more explicit than the following definitions.
The condition $\kappa$ defines two points $x_n$ and $y_n$ on the boundary of
$\cyl(nA, h(n))$ (see section \ref{chapitre6subsubsec:lowerboundpgd}). As in section \ref{chapitre6subsubsec:lowerboundpgd}, we denote
by $L(n,\wt)$ the distance between $x_n$ and $y_n$, and we have
$$ L(n,\wt) \,=\, \frac{nl(A)}{cos(\widetilde{\theta}-\theta)}\,. $$ We
define
$$ z_i \,=\, \left(\zeta'(n) + (i-1) L(n,\wt)\right) \vd \,,$$
for $i=1,..., \N$, where
$$ \N \,=\, \left\lfloor \frac{N-2\zeta'(n)}{L(n,\wt)} \right\rfloor \,.$$
Of course we consider only $N$ large enough to have $\N \geq 2$. For
$i=1,...,\N$, we denote by $\widetilde{B_i}$ the image of $\cyl
(nA, h(n))$ by the translation of vector $\overrightarrow{x_n z_i}$. For $N$
sufficiently large, thanks to condition (\ref{chapitre6cond2}), we know that
$\widetilde{B_i} \subset \cyl''(N)$ for all $i$. We can translate $\widetilde{B_i}$ again by a vector of
norm strictly smaller than $1$ to obtain an integer translate of $\cyl(nA,
h(n))$ (i.e., a translate by a vector whose coordinates are in $\ZZ^2$) that we will call $B_i$. Now we want to glue together
cutsets of boundary condition $\kappa$ in the different $B_i$. We define:
$$ \E_1(n,N,\kappa) \,=\, \left(\bigcup_{i=1}^\N \mathcal{V} (z_i, \zeta) \right)
\, \bigcap \,\cyl''(N)  \,,$$
where $\zeta$ is still a fixed constant bigger than $4$, and:
$$ \E_2(n,N,\kappa) \,=\, \mathcal{V} \left( [0,\zeta'(n)
  \vd ] \cup [z_\N, N \vd ]
  ,\zeta \right)\, \bigcap \,\cyl''(N) \,.$$
Let $E_1(n,N,\kappa)$ (respectively $E_2(n,N,\kappa)$) be the set of the edges included in
$\E_1(n,N,\kappa)$ (respectively $\E_2(n,N,\kappa)$).
Then, still by gluing cutsets together, we obtain:
\begin{equation}
\label{chapitre6lien_phiF_tau_2}
\tau(\cyl''(N), \vb ) \,\leq\, \sum_{i=1}^\N
\phi^{\kappa}(B_i, \va) + V(E_1(n,N,\kappa) \cup E_2(n,N,\kappa)) \,.
\end{equation}
This equation (\ref{chapitre6lien_phiF_tau_2}) is equation (24) in \cite{RossignolTheret09}.
On one hand, there exists a constant $C_6$ (independent of $\kappa$) such that:
$$ \card (E_1(n,N,\kappa) \cup E_2(n,N,\kappa)) \,\leq\, C_6 \left( \N + \zeta'(n) +
  L(n,\wt) \right) \,.$$ 
On the other hand, the variables $(\phi^{\kappa}(B_i))_{i=1,...,\N}$ are
identically distributed, with the same law as $\phi_{n}^\kappa$ (because we
only consider integer translates). Then for all $\kappa \in D_n$, for all $\widetilde{\lambda}\geq
\varepsilon>0$, for all large $N$, we have by the FKG inequality
\begin{align*}
\PP \Bigg[&\tau(\cyl''(N)  , \vb ) \leq
  \left(\widetilde{\lambda} -\frac{1}{\sqrt{N}}  \right) N \Bigg] \\
& \,\geq\, \PP \left[ \phi_n^\kappa \leq \left(
    \widetilde{\lambda}-\varepsilon \right) \frac{N}{\N} \right]^{\N}
\times \PP \left[V(E_1(n,\kappa)\cup E_2(n,\kappa)) \leq
  \frac{\varepsilon}{2} N\right] \\
&\,\geq\, \PP \left[ \frac{\phi_n^\kappa}{nl(A)} \leq
  \frac{\widetilde{\lambda} - \eps}{ \cos (\wt - \theta)} \right]
^{\N}\times \PP
\left[ t(e) \leq \frac{\varepsilon N}{2 C_6 (\N + \zeta'(n) + L(n,\wt))}
\right]^{C_6 (\N + \zeta'(n) + L(n,\wt))}\\ 
& \,\geq\, \PP \left[ \frac{\phi_n^\kappa}{nl(A)} \leq
  \frac{\widetilde{\lambda} -\eps}{ \cos (\wt - \theta)} \right]^{\N} \times \PP
\left[ t(e) \leq \frac{\varepsilon l(A) n}{4 C_6 }
\right]^{C_6 (\N + \zeta'(n) + L(n,\wt))} \,.
\end{align*}
We take the logarithm of the previous inequality, divide it by $-N$, and
send $N$ to infinity to obtain that:
$$ \I_{\wt} (\widetilde{\lambda}) \,\leq\, \frac{-1}{L(n,\wt)} \log \PP
\left[ \frac{\phi_n^\kappa}{nl(A)} \leq
  \frac{\widetilde{\lambda} -\eps }{ \cos (\wt - \theta)} \right]  -
\frac{C_6}{L(n,\wt)} \log \PP \left[ t(e) \leq \frac{\varepsilon l(A)
  n}{4C_6} \right] \,. $$
For $n$ large enough,
$$\PP\left[ t(e) \leq \frac{\varepsilon l(A)
  n}{4C_6} \right] \geq \frac{1}{2} \,,$$
and thus,
$$\frac{1}{nl(A)} \log \PP \left[ \frac{\phi_n^\kappa}{nl(A)} \leq
  \frac{\widetilde{\lambda} -\varepsilon}{\cos (\wt - \theta)} \right]
\,\leq\, - \frac{1}{\cos (\wt - \theta)} \I_{\wt} (\widetilde{\lambda}) +
\frac{K}{n}\,,$$
where $K = C_6 \log 2 / l(A)$. We set $\lambda = (\widetilde{\lambda} -
\varepsilon)/ \cos(\wt - \theta)$ (so $\lambda \geq 0$), and  let
$\varepsilon$ go to zero to conclude that for all $\lambda \geq 0$ and $\kappa \in D_n$, 
\begin{equation}
\label{chapitre6dec}
\PP \left[\frac{\phi_n^\kappa}{nl(A)} \leq \lambda \right] \,\leq\, \exp
- \left[ n l(A) \left( \frac{1}{\cos(\wt -\theta)} \I_{\wt} (\lambda
    \cos(\wt -  \theta) ^+) + \frac{K}{n} \right) \right]\,.
\end{equation}
We come back now to the study of $\phi_n$ itself. We have seen that $\phi_n
= \inf_{\kappa \in D_n} \phi_n^\kappa$. We also noticed that
$\phi_n^\kappa$ takes only a finite number of values when $\kappa \in D_n$,
thus one may restrict ourselves to a finite subset $\tilde{D}_n$ of $D_n$
such that $\card (\tilde{D}_n) \leq C_4 h(n)^2$. Therefore,
\begin{align*}
\PP [\phi_n\leq \lambda n l(A)] & \,=\, \PP \left[ \exists \kappa \in
  \widetilde{D}_n'
\,|\, \phi_n^\kappa \leq \lambda n l(A)\right]\\
& \,\leq\, \sum_{\kappa \in \widetilde{D}_n'} \PP [\phi_n^\kappa \leq \lambda n
l(A)]\\
& \,\leq\, C_4 h(n)^2 \times \max_{\kappa \in D_n} \PP [\phi_n^\kappa \leq
\lambda n l(A)]\\
& \,\leq\, C_4 h(n)^2 \exp \left[-n l(A) \left( \frac{K}{n} + \inf_{\wt
      \in \D_n} \frac{1}{\cos (\wt - \theta)} \I_{\wt} (\lambda \cos (\wt -
    \theta )^+)  \right)\right]\,.
\end{align*}
If we suppose that $\lim_{n\rightarrow \infty} \log (h(n)) / n =0$, we obtain that:
\begin{equation}
\label{chapitre6interm}
\limsup_{n\rightarrow \infty} \frac{1}{nl(A)} \log \PP [\phi_n \leq
\lambda nl(A) ] \,\leq\, - \liminf_{n\rightarrow \infty} \inf_{\wt \in \D_n}
\frac{1}{\cos(\wt - \theta)} \I_{\wt} (\lambda \cos(\wt - \theta)^+)\,,
\end{equation}
and
\begin{equation}
\label{chapitre6intermbis}
\liminf_{n\rightarrow \infty} \frac{1}{nl(A)} \log \PP [\phi_n \leq
\lambda nl(A) ] \,\leq\, - \limsup_{n\rightarrow \infty} \inf_{\wt \in \D_n}
\frac{1}{\cos(\wt - \theta)} \I_{\wt} (\lambda \cos(\wt - \theta)^+)\,,
\end{equation}
We can now apply Lemmas 4.1 and 4.2 in \cite{RossignolTheret09} with $f=g_\lambda$, that we know to be
l.s.c. thanks to Lemma \ref{chapitre6lem1}, to obtain that
$$ \liminf_{n\rightarrow \infty} \inf_{\wt \in \D_n} g_{\lambda} (\wt)
\,\geq\, \inf_{\wt \in \ad(\overline{\D}) }g_{\lambda}(\wt) \quad \textrm{and} \quad
\limsup_{n\rightarrow \infty} \inf_{\wt \in \D_n} \,\geq\, \inf_{\wt \in
  \ad(\underline{\D}) } g_{\lambda}(\wt)\,.   $$
So (\ref{chapitre6interm}) and (\ref{chapitre6intermbis}) lead to (\ref{chapitre6eq2}) and (\ref{chapitre6eq2bis}).


\subsection{Positivity and proof of Theorem~\ref{thmLDphi}}
\label{secpositivite}

We need some extra hypotheses: from now on, we suppose that $\mathbf{(H1)}$,
$\mathbf{(H2)}$, $\mathbf{(F1)}$, $\mathbf{(F2)}$ and $\mathbf{(FH1)}$
hold. We define 
$$\eta_{\theta,h}\,=\, \inf_{\wt \in \overline{\D}} \frac{\nu_{\wt}}{\cos (\wt - \theta)}
\,=\,  \inf_{\wt \in \underline{\D}} \frac{\nu_{\wt}}{\cos (\wt - \theta)}\,. $$
We know that under these hypotheses,
$$ \lim_{n\rightarrow \infty} \frac{\phi_n}{nl(A)} \,=\, \eta_{\theta,h} \,>\,0
\quad \textrm{in }L^1 \,.$$
In fact, the convergence of $\phi_n / (nl(A))$ in $L^1$ is stated in
Corollary \ref{chapitre6corollaire} under the stronger assumption
$\mathbf{(H3)}$. However, the methods used in \cite{RossignolTheret08b} to
prove this $L^1$-convergence do not use the condition
$\mathbf{(H3)}$ itself, and can be performed under the assumption
$\mathbf{(FH1)}$ instead.

We define
$$\wK(\lambda) \,=\,  \inf_{\wt \in \ad(\overline{\D})}
\frac{1}{\cos (\wt - \theta)} \I_{\wt} (\lambda \cos(\wt
-\theta)^+)\,=\, \inf_{\wt \in \ad(\overline{\D})} g_{\lambda} (\wt)\,.   $$
We want to prove that $\wK>0$ on $[0,\eta_{\theta,h}[$. Indeed,
we know that
$$ \limsup_{n\rightarrow \infty} \frac{1}{nl(A)} \log \PP \left[ \phi_n
  \leq \lambda n l(A) \right] \,\leq\, - \wK \,, $$
thus proving that $\wK >0$ on $[0,\eta_{\theta,h}[$ is
equivalent to proving that the lower large
deviations of $\phi_n/(nl(A))$ are (at least) of surface order.

We remember all the properties of $\I_{\wt}$ we know (see Theorem
\ref{propI}). Let us define
$$\delta_{\theta, h} \,=\, \delta \times  \inf_{\wt \in \D} \frac{|\cos \wt| +
|\sin \wt|}{\cos (\wt - \theta)}\,, $$
the infimum of the values that $\phi_n/(nl(A))$ can take asymptotically.
Then obviously $\wK$ is infinite on $[0, \delta_{\theta, h}[$
and finite on $]\delta_{\theta, h}, + \infty[$ (its
behaviour at $\delta_{\theta, h}$ will be study in section \ref{chapitre6propK}). It is also obvious that
$\widetilde{K}$ is null on
$[\eta_{\theta, h}, +\infty[$. Since for all $\wt$, $\lambda \rightarrow
\I_{\wt}(\lambda^+)$ is non increasing, so is $\wK$ on
$\RR^+$. We state the following result:
\begin{lem}
\label{chapitre6strict}
The function $\wK$ is strictly decreasing on
$[\delta_{\theta,h},\eta_{\theta, h}]$, i.e.,
$$ \forall \lambda \in ]\delta_{\theta,h},\eta_{\theta, h}] \,,\,\,
\forall \eps >0 \qquad \wK (\lambda) \,<\, \wK
(\lambda - \eps) \,.$$
\end{lem}
We immediately notice that this lemma implies the
positivity of $\wK$ on $[0,\eta_{\theta,h}[$, and thus Theorem~\ref{thmLDphi}
through inequality (\ref{chapitre6interm}).

\begin{dem}
Thanks to Lemma \ref{chapitre6lem1}, we know that for every fixed $\lambda$,
$\wt\mapsto g_\lambda(\wt)$ is l.s.c. Since $\ad(\overline{\D})$ is compact,
$\inf_{\wt \in \ad(\overline{\D})} g_{\lambda} (\wt)$ is reached at some
$\wt_\lambda\in\ad(\overline{\D})$. Notice also that for every fixed $\wt$,
$\lambda\mapsto g_\lambda(\wt)$ is strictly
decreasing (in the same meaning as in Lemma \ref{chapitre6strict}) on
the interval:
$$\left[\delta\frac{|\cos \wt | + |\sin \wt |}{\cos (\wt - \theta)},\frac{\nu_{\wt} }{ \cos(\wt -
\theta)}\right]\;.$$ 
We consider $\lambda \in ]\delta_{\theta,h},\eta_{\theta,
  h}]$. Thus $\wK(\lambda) <\infty$, so we can suppose that
$\wK(\lambda - \eps) <\infty$ otherwise the result is
obvious. The condition $\wK(\lambda - \eps) <\infty$ is
equivalent by definition of $\wt_{\lambda -\eps}$ to $g_{\lambda - \eps}
(\wt_{\lambda - \eps}) <\infty$, which implies that
$$ \lambda \,>\, \delta \frac{|\cos \wt_{\lambda - \eps}| + |\sin
  \wt_{\lambda - \eps}|}{\cos (\wt_{\lambda - \eps} - \theta)}\,. $$
We deduce that 
$$\lambda \in \left] \delta\frac{|\cos \wt_{\lambda - \eps} | + |\sin \wt_{\lambda - \eps} |}{\cos (\wt_{\lambda - \eps} - \theta)},\frac{\nu_{\wt_{\lambda - \eps}} }{ \cos(\wt_{\lambda - \eps} -\theta)}\right]\,, $$
thus
$$ g_{\lambda} (\wt_{\lambda - \eps}) \,<\, g_{\lambda-\eps} (\wt_{\lambda
  - \eps})\,. $$
We obtain
$$ \wK(\lambda)\,=\,\inf_{\wt\in\ad(\overline{\D})}g_\lambda(\wt)\,\leq\,
g_\lambda(\wt_{\lambda-\eps})\,<\,g_{\lambda-\eps}
(\wt_{\lambda-\eps})\,=\,\inf_{\wt\in\ad(\overline{\D})}g_{\lambda-\eps}(\wt)\,=\,
\wK (\lambda - \eps)\;.$$
\end{dem}


\subsection{Discussion}
\label{secdiscussion}

Combining the results of the two previous sections, we obtain that for all
$\lambda \geq 0$, if we define
$$ \square_n \,=\, \frac{1}{nl(A)} \PP [\phi_n \leq \lambda n l(A)]\,, $$
we have
$$ \left\{\begin{array}{l} -\begin{displaystyle} \inf_{\wt \in
        \ad(\underline{\D})}  \end{displaystyle}\frac{ \I_{\wt} (\lambda \cos (\wt - \theta) ^-)}{\cos
      (\wt - \theta)} \,\leq\,
   \begin{displaystyle}  \liminf_{n\rightarrow \infty} \end{displaystyle} \square_n \,\leq\, -\inf_{\wt \in
      \ad(\underline{\D})} \frac{\I_{\wt} (\lambda \cos
    (\wt - \theta) ^+)}{\cos (\wt - \theta)} \,, \\
 - \begin{displaystyle}\inf_{\wt \in \ad(\overline{\D})} \end{displaystyle} \frac{ \I_{\wt} (\lambda \cos (\wt - \theta) ^-)}{\cos
      (\wt - \theta)} \,\leq\,
    \limsup_{n\rightarrow \infty} \square_n \,\leq\, - \begin{displaystyle}\inf_{\wt \in
      \ad(\overline{\D})} \end{displaystyle} \frac{\I_{\wt} (\lambda \cos
    (\wt - \theta) ^+)}{\cos (\wt - \theta)} \, .
  \end{array} \right. $$
In fact, we will prove in section \ref{chapitre6propK} that for all
$$ \lambda \,\neq\, \delta \inf_{\wt \in \D} \frac{|\cos \wt|+|\sin \wt|}{\cos
(\wt -\theta)} \,=\, \delta_{\theta,h}\,, $$
we have
\begin{equation}
\label{eqegaliteinf}
 \inf_{\wt \in \D}\frac{ \I_{\wt} (\lambda \cos (\wt -
\theta) ^+)}{\cos (\wt - \theta)} \,=\,  \inf_{\wt \in \D}\frac{ \I_{\wt}
(\lambda \cos (\wt - \theta) ^-)}{\cos (\wt - \theta)}\,,
\end{equation}
for $\D$ equal to $\ad(\underline{\D})$ or $\ad(\overline{\D})$ (see the
proof of the continuity of $\wK$, Lemma \ref{continuiteK}). It implies
that for all $\lambda \neq \delta_{\theta,h}$,
$$ \liminf_{n\rightarrow \infty} \square_n  \,=\,  -\inf_{\wt \in \ad(\underline{\D})} \frac{1}{\cos
      (\wt - \theta)} \I_{\wt} (\lambda \cos (\wt - \theta) ^+)\,=\,
    -\inf_{\wt \in \ad(\underline{\D})} \frac{1}{\cos  (\wt - \theta)} \I_{\wt} (\lambda
    \cos (\wt - \theta) ^-) $$
and
$$ \limsup_{n\rightarrow \infty} \square_n  \,=\,  -\inf_{\wt \in \ad(\overline{\D})} \frac{1}{\cos
      (\wt - \theta)} \I_{\wt} (\lambda \cos (\wt - \theta) ^+)\,=\,
    -\inf_{\wt \in \ad(\overline{\D})} \frac{1}{\cos  (\wt - \theta)} \I_{\wt} (\lambda
    \cos (\wt - \theta) ^-) \,,$$
thus under condition $\mathbf{(FH2)}$ we obtain that for all $\lambda \neq
\delta_{\theta,h}$, we know that $\lim_{n\rightarrow \infty} \square_n$
exists and
$$ \lim_{n\rightarrow \infty} \frac{1}{nl(A)}\PP [\phi_n \leq \lambda n
l(A)] \,=\, \wK (\lambda) \,.  $$
Notice that if $\I_{\wt_0} (\delta(|\cos \wt_0|+|\sin \wt_0|) ^+)
<\infty$ for some $\wt_0 \in \D$ such that $\delta_{\theta,h} = \delta
(|\cos \wt_0|+|\sin \wt_0|)/\cos (\wt_0 -\theta) $, then
$$ -\infty \,=\, -\inf_{\wt \in \D} \frac{1}{\cos
      (\wt - \theta)} \I_{\wt} (\delta_{\theta,h} \cos (\wt - \theta) ^-) \,<\, 
-\inf_{\wt \in \D} \frac{1}{\cos (\wt - \theta)} \I_{\wt}
(\delta_{\theta,h} \cos (\wt - \theta) ^+)  \,,$$
thus we have no hope to prove equation (\ref{eqegaliteinf}) at the point
$\lambda = \delta_{\theta,h}$.

\section{Large deviation principle}
\label{chapitre6subsec:LDP}

From now on, we suppose that hypotheses $\mathbf{(H1)}$, $\mathbf{(H2)}$,
$\mathbf{(F1)}$, $\mathbf{(F2)}$, $\mathbf{(FH1)}$ and $\mathbf{(FH2)}$
hold (we recall that $\mathbf{(H4)}$ implies $\mathbf{(FH1)}$ and
$\mathbf{(FH2)}$). By definition, we have
\begin{align*}
\forall \lambda \in \RR^+ ,\, \wK(\lambda) & \,=\,  \inf_{\wt
  \in \ad(\overline{\D})}
\frac{1}{\cos (\wt - \theta)} \I_{\wt} (\lambda \cos(\wt  -\theta)^+)\\
& \,=\,\inf_{\wt \in \ad(\underline{\D})} \frac{1}{\cos (\wt - \theta)}
\I_{\wt} (\lambda \cos(\wt  -\theta)^+)\,.
\end{align*}
We define the rate function $\K: \RR^+ \rightarrow \RR^+ \cup
\{+\infty\}$ by:
$$ \K(\lambda) \,=\,\left\{ \begin{array}{ll} \inf_{\wt \in \ad(\overline{\D})}
    \frac{1}{\cos (\wt
    - \theta)} \I_{\wt} (\lambda \cos(\wt  -\theta)^+) \,=\, \widetilde{K}(\lambda)  & \,\, \mbox{if} \,\,
  \lambda \leq \eta_{\theta,h}\,,\\
+\infty & \,\, \mbox{if} \,\,  \lambda > \eta_{\theta,h}  \,.\end{array}
\right.$$


\subsection{Properties of $\K$}
\label{chapitre6propK}

We can deduce a
lot of properties of $\K$ from the properties of $\wK$ stated
in section \ref{secpositivite}: $\K$ is infinite on $[0, \delta_{\theta,
  h}[$ (if $\delta_{\theta,h}>0$) and on $]\eta_{\theta, h}, + \infty[$,
finite on $]\delta_{\theta, h}, \eta_{\theta, h}]$ (if $\eta_{\theta, h} >
\delta_{\theta, h}$) and strictly decreasing on $[\delta_{\theta,h}, \eta_{\theta, h}]$ in
the sense of Lemma \ref{chapitre6strict}. We only have to prove that $\K$ is a
good rate function, and that it is continuous on $[0,\eta_{\theta,h}]$
except possibly at $\delta_{\theta,h}$ where it may be only right
continuous. We first state that $\K$ is a good rate function:
\begin{lem}
\label{chapitre6Ksci}
The function $\K$ is lower semi-continuous and coercive on $\RR^+$, i.e.,
for all $t \geq 0$, the set $ \{  \lambda \,|\, \K(\lambda) \leq t \}$ is
compact.
\end{lem}
We will use this property to prove that $\wK$ is right continuous.
\begin{dem}
In fact it is sufficient to prove that for all $t\geq 0$, the set $\{
\lambda \,|\, \wK (\lambda) \leq t \}$ is closed, because we
know that
$$ \forall t\geq 0 \qquad  \{  \lambda \,|\, \K(\lambda) \leq t \} \,=\, \{
\lambda \,|\, \wK (\lambda) \leq t \} \cap [0,\eta(\theta,
h)] \,. $$
Let $(\lambda_n)_{n\geq 0}$ be a sequence of $\{
\lambda \,|\, \wK (\lambda) \leq t \}$, converging towards some
$\lambda_0$. For each fixed $\lambda$ in $\RR^+$, since the function
$g_{\lambda}$ is lower semi-continuous and $\ad(\overline{\D})$ is compact,
there exists 
$\wt_{\lambda}$ such that
$$  \wK(\lambda) \,=\, g_{\lambda} (\wt_{\lambda}) \,.$$
The sequence $(\wt_{\lambda_n})_{n\geq 0}$ takes values in the compact
$\ad(\overline{\D})$, so up to extracting a subsequence, we can suppose that it converges towards a limit  $\wt_0 \in \ad(\overline{\D})$. For all positive
$\eps$, for all large $n$ we
have $\lambda_n \leq \lambda_0 + \eps$, and so, since $\I_{\wt}$ is non
increasing for all $\wt$, we obtain for all large $n$ that
$$ g_{(\lambda_0 + \eps)}(\wt_{\lambda_n}) \,\leq\, g_{\lambda_n}
(\wt_{\lambda_n})\,\leq\, t  \,.$$
Since $g_{(\lambda_0 +\eps)}$ is l.s.c. and a subsequence
$(\wt_{\psi(n)})_{n\geq 0}$ of
$(\wt_{\lambda_n})_{n\geq 0}$ converges towards $\wt_0$, we obtain:
$$ g_{(\lambda_0 + \eps) } (\wt_0) \,\leq\, \liminf_{n\rightarrow \infty}
g_{(\lambda_0+\eps )} (\wt_{\psi(n)}) \,\leq\, t\,. $$
This inequality is satisfied for all positive $\eps$, and $\wt_0 \in \ad(\overline{\D})$, so
$$ \wK(\lambda_0)\,\leq\, g_{\lambda_0}(\wt_0) \,=\, \lim_{\eps
  \rightarrow   0 \,,\, \eps >0} g_{(\lambda_0+\eps)} (\wt_0)\,\leq\, t\,.  $$
This ends the proof of Lemma \ref{chapitre6Ksci}.
\end{dem}

We now study the continuity of $\K$:
\begin{lem}
\label{continuiteK}
The function $\wK$ is continuous on $\mathbb{R}^+$, except possibly at
$\delta_{\theta,h}$ where it may be only right continuous.
\end{lem}
The proof of the continuity of $\wK$ is quite long and technical, and this
property of $\wK$ is not needed to prove the large deviation
principle. However, as explained in section \ref{secdiscussion} and below,
the continuity of $\wK$ is a natural question to ask, so it seems to us
important to give an answer to it.
\begin{dem}
We define
$$ \wK (\lambda ^+) \,=\, \lim_{\eps \rightarrow 0,\, \eps >0} \wK (\lambda
+ \eps) \quad \textrm{and} \quad  \wK (\lambda ^-) \,=\, \lim_{\eps
  \rightarrow 0,\, \eps >0} \wK (\lambda - \eps)\,. $$
First of all, we prove that $\wK$ is right continuous, i.e., $\wK(\lambda
^+)=\wK (\lambda)$ for all $\lambda \in \RR^+$. We have for all $\lambda
\geq 0$
\begin{align*}
\wK (\lambda ^+ ) & \,=\, \lim_{\eps \rightarrow 0, \, \eps >0} \inf_{\wt
  \in \ad (\overline{\D})} \frac{1}{\cos (\wt - \theta)} \I_{\wt} ((\lambda
+ \eps) \cos (\wt - \theta) ^+)\\
& \,\leq\,  \lim_{\eps \rightarrow 0, \, \eps >0} \inf_{\wt
  \in \ad (\overline{\D})} \I_{\wt} (\lambda \cos (\wt - \theta) ^+)\\
& \,\leq\, \wK (\lambda)\,,
\end{align*}
and since $\wK$ is lower semi-continuous, we know that
$$ \wK (\lambda ^+) \,\geq\, \wK(\lambda) \,, $$
thus
$$ \forall \lambda \geq 0\quad \wK(\lambda ^+) \,=\, \wK (\lambda)\,. $$
On the other hand, for all $\lambda \geq 0$ we have
\begin{align*}
\wK(\lambda ^-) & \,=\, \lim_{\eps \rightarrow 0,\, \eps >0} \inf_{\wt \in
  \ad(\overline{\D})}  \frac{1}{\cos (\wt - \theta)} \I_{\wt} ((\lambda
- \eps) \cos (\wt - \theta) ^+)\\
& \,=\, \inf_{\wt \in \ad(\overline{\D})} \lim_{\eps \rightarrow 0,\, \eps
  >0}  \frac{1}{\cos (\wt - \theta)} \I_{\wt} ((\lambda - \eps) \cos (\wt - \theta) ^+)\\
& \,=\, \inf_{\wt \in \ad(\overline{\D})} \frac{1}{\cos (\wt - \theta)}
\I_{\wt} (\lambda \cos (\wt - \theta) ^-)\\
& \,\geq\, \wK(\lambda) \,,
\end{align*}
since the limit in $\eps$ appearing in these equations is a decreasing
limit. Thus $\wK$ is continuous at $\lambda$ if and only if
$$  \inf_{\wt \in \ad(\overline{\D})} \frac{1}{\cos (\wt - \theta)}
\I_{\wt} (\lambda \cos (\wt - \theta) ^-) \,=\,  \inf_{\wt \in
  \ad(\overline{\D})} \frac{1}{\cos (\wt - \theta)} \I_{\wt} (\lambda \cos
(\wt - \theta) ^+) \,. $$
We will thus prove that it is true for all $\lambda \neq
\delta_{\theta,h}$. Notice that the proof we propose can be performed with
$\ad(\underline{\D})$ instead of $\ad(\overline{\D})$, and thus the
continuity of $\wK$ is linked with the existence of the limit
$$ \lim_{n\rightarrow \infty} \frac{1}{nl(A)}\PP [\phi_n \leq \lambda n
l(A)] $$
as explained in section \ref{secdiscussion}. We need two intermediate
lemmas to prove Lemma \ref{continuiteK}. The first one is the following:
\begin{lem}
\label{lemconvexiteLambda}
Let $\Lambda$ be the function defined on $\RR^+ \times
\RR^2\smallsetminus{(0,0)}$ by
$$ \Lambda (\lambda, \vec{v}) \,=\,\|\vec{v}\|_2  \I_{\theta(\vec{v})} \left( \lambda
  (|\cos(\theta(\vec{v}))| + |\sin (\theta(\vec{v}))| )^+ \right)
 $$
where $\theta(\vec{v}) \in [0,2\pi[$ satisfies $\vec{v} = \|\vec{v}\|_2 (\cos
(\theta(\vec{v})), \sin (\theta(\vec{v})))$ and $\|\vec{v}\|_2$ is the
Euclidean norm of $\vec{v}$. Then for all vectors $\vec{u}$ and $\vec{v}$
in $(\RR^+)^2\smallsetminus \{(0,0)\}$, we have
\begin{equation}
\label{eqconvexite}
\forall \lambda \in \RR^+ \qquad \Lambda (\lambda, \vec{u} + \vec{v})
\,\leq\, \Lambda (\lambda, \vec{u}) + \Lambda (\lambda, \vec{v})\,.
\end{equation}
\end{lem}
\begin{dem}
Lemma \ref{lemconvexiteLambda} is a simple consequence of Lemma
\ref{chapitre6trigI}. We consider $\vec{u}$ and $\vec{v}$ in
$(\RR^+)^2\smallsetminus \{(0,0)\}$, and define $\vec{w} = \vec{u} +
\vec{v}$. We use the notations $\theta(\vec{u}) = \wt_a$, $\theta(\vec{v})
= \wt_b$ and $\theta(\vec{w}) = \wt_c$. We consider the triangle $(abc)$ of
side $[bc]$ (resp. $[ab]$, $[ac]$) orthogonal to $\vec{u}$
(resp. $\vec{w}$, $\vec{v}$) and of length $\|\vec{u}\|_2$
(resp. $\|\vec{w}\|_2$, $\|\vec{v}\|_2$). It is indeed a triangle since
$\vec{w} = \vec{u} + \vec{v}$. Moreover, since $\vec{u}$ and $\vec{v}$ are
in $(\RR^+)^2\smallsetminus \{(0,0)\}$, we know that the angles
$\widehat{cab}$ and $\widehat{abc}$ have values strictly smaller than
$\pi/2$. We consider $\lambda \in \RR^+$, and
$$ \lambda ' \,=\, \lambda\, l(ab) (\cos \wt_c + \sin \wt_c) \,=\, \lambda\,
\|\vec{w}\|_1 \,. $$
We can apply Lemma \ref{chapitre6trigI} in the triangle
$(abc)$ to obtain for all $\alpha \in [0,1]$
\begin{align*}
\Lambda (\lambda, \vec{u} + \vec{v}) & \,=\, l(ab) \, \I_{\wt_c} \left(
  \frac{\lambda'}{l(ab) }^+  \right)\\
& \,\leq\, l(ac)\, \I_{\wt_b} \left(\alpha \frac{\lambda '}{l(ac)} ^+ \right)
+  l(bc)\, \I_{\wt_a} \left((1-\alpha) \frac{\lambda '}{l(bc)} ^+ \right)\\
& \,\leq\, \Lambda \left( \lambda \frac{\alpha \|\vec{w}\|_1}{\|\vec{v}\||1} , \vec{v}\right) + \Lambda \left( \lambda \frac{(1-\alpha) \|\vec{w}\|_1}{\|\vec{u}\||1} ,\vec{u}\right)\,.
\end{align*}
Since $\vec{u}$ and
$\vec{v}$ are both in $(\RR^+)^2$, we know that $\|\vec{w}\|_1 =
\|\vec{u}\|_1 + \|\vec{v}\|_1$, thus we can choose $\alpha = \|\vec{v}\|_1
/ \|\vec{w}\|_1 \in [0,1]$ and $1-\alpha = \|\vec{u}\|_1 /
\|\vec{w}\|_1$. This ends the proof of Lemma \ref{lemconvexiteLambda}.
\end{dem}

Lemma \ref{lemconvexiteLambda} states a property of convexity for the
function $\Lambda$. To deduce from it a property of continuity, we need to
investigate when $\Lambda$, thus $\I_{\wt}$, is finite. It is well known
(see \cite{RossignolTheret08b}) that $\I_{\wt}(\lambda ^+) $ is infinite if
$\lambda \in [0, \delta
(|\cos \wt| + |\sin \wt|)[$ (if $\delta>0$) and finite if $\lambda \in]\delta (|\cos
\wt| + |\sin \wt|), +\infty[ $. The only point to study is the
behaviour of $\I_{\wt}(\delta (|\cos \wt| + |\sin \wt|)^+)$. This is the
purpose of the following Lemma (which could be stated in dimension $d\geq
2$ in fact):
\begin{lem}
\label{lemfini}
Whatever the value of $\delta$, we have
\begin{align*}
 \forall \wt \,,\,\, \I_{\wt} (\delta (|\cos \wt| + |\sin \wt|)^+) \,<\,
\infty \quad \iff \quad& \exists  \wt \,,\,\, \I_{\wt} (\delta (|\cos \wt| + |\sin \wt|)^+) \,<\,
\infty \\& \quad \iff \quad \PP(t(e) = \delta) \,>\, 0\,.
\end{align*}
\end{lem}
\begin{dem}
First, let us prove that
\begin{equation}
\label{eqfini1}
\PP (t(e) = \delta) \,>\,0 \quad \Longrightarrow  \quad \forall \wt\,,\,\, \I_{\wt} (\delta (|\cos
\wt| + |\sin \wt|)^+ ) \,<\, \infty\,.
\end{equation}
Let $\widetilde{A}$ be a line segment orthogonal to $(\cos \wt, \sin \wt)$
for some fixed $\wt$, and $h$ a height function satisfying
$\lim_{n\rightarrow \infty} h(n) = +\infty$. We know (see Lemma 4.1 in \cite{RossignolTheret08b}) that the minimal
number of edges $\N(n\widetilde{A},h(n))$ of a cutset that separates
$(n\widetilde{A})_1^{h(n),1/2,\wt}$ from
$(n\widetilde{A})_2^{h(n),1/2,\wt}$  satisfies
\begin{equation}
\label{eqN}
 \left| \frac{\N(n\widetilde{A},h(n))}{nl(\widetilde{A})} - (|\cos \wt| +
  |\sin \wt|) \right| \,\leq\, \frac{2}{nl(\widetilde{A})} \,.
\end{equation}
Let $\eps>0$. Let $E_{\min}(n)$ be a cutset of minimal number of edges. For $n$
large enough, we know that $1/\sqrt{n} \leq \eps /2$ and
$\N(n\widetilde{A},h(n))/(nl(\widetilde{A})) \leq (|\cos \wt| + |\sin \wt|)
+ \eps /(2\delta)$, and we obtain
\begin{align*}
\PP \Bigg( \frac{\tau(n\widetilde{A},h(n))}{nl(\widetilde{A})} \leq & \delta
  (|\cos \wt| + |\sin \wt| ) + \eps - \frac{1}{\sqrt{n}} \Bigg)\\ & \,\geq
\, \PP \left( \frac{V(E_{\min}(n))}{nl(\widetilde{A})} \leq  \delta (|\cos
  \wt| + |\sin \wt| ) + \frac{\eps}{2}  \right)\\
& \,\geq\, \PP \left(\forall e\in E_{\min}(n)\,,\,\, t(e) =\delta  \right)\\
&\,\geq\, \PP(t(e) = \delta)^{\N(n\widetilde{A},h(n))}\,.
\end{align*}
Thus for all $\eps >0$ we have
$$ \I_{\wt} (\delta (|\cos \wt| + |\sin \wt|) + \eps) \,\leq\, - (|\cos
\wt| + |\sin \wt|) \log \PP(t(e) = \delta) \,, $$
and we conclude that
$$ \I_{\wt} (\delta (|\cos \wt| + |\sin \wt|)^+) \,\leq\,  - (|\cos
\wt| + |\sin \wt|) \log \PP(t(e) = \delta)\,. $$
This implies (\ref{eqfini1}). We now prove that if for all distribution
function $G$ such that $\inf \{ x\,|\, G(x)>0 \} =0$, we have
$$G(0) \,=\, 0 \quad \Longrightarrow \quad \forall
\wt\,,\,\,\I^{[G]}_{\wt} (0^+) \,=\, +\infty \,,$$
where the exponent $[G]$ stress the dependence of $\I$ in $G$,
then for all distribution function $F$ on $\RR^+$, if $\delta=\inf
\{ x\,|\, F(x)>0 \}$, we obtain:
$$ F(\delta) \,=\, 0 \quad \Longrightarrow \quad
\forall \wt\,,\,\,\I^{[F]}_{\wt} (\delta (|\cos \wt| + |\sin \wt|) ^+) \,=\,
+\infty \,.$$
Let $F$ be a distribution function on $\RR^+$, $\delta =\inf \{ x\,|\, F(x)>0 \} $, and
$(t(e))$ the family of capacities on the edges of distribution function
$F$. Let $t'(e) = t(e) - \delta \geq 0$ for all $e$, $t'(e)$ has distribution
function $G$ such that $\inf \{ x\,|\, G(x)>0 \} =0$, and $G(0) =
F(\delta)$. We denote by $\tau$ (resp. $\tau'$) the maximal flows
corresponding to the capacities $(t(e))$ (resp. $(t'(e))$). Then obviously
$$ \tau(n\widetilde{A}, h(n)) \,\geq\, \tau'(n\widetilde{A}, h(n)) +
\delta  \N (n\widetilde{A}, h(n)) \,,$$
thus for $n$ large enough, thanks to (\ref{eqN}), we have
$$ \PP \left( \frac{\tau(n\widetilde{A}, h(n))}{nl(\widetilde{A})} \leq \delta
  (|\cos \wt| + |\sin \wt|) + \eps - \frac{1}{\sqrt{n}}  \right) \,\leq\,
\PP \left(\frac{\tau'(n\widetilde{A}, h(n))}{nl(\widetilde{A})} \,\leq
  \frac{\eps}{4} \right)\,. $$
Thus
$$ \I_{\wt}^{[F]} (\delta (|\cos \wt| + |\sin \wt|)^+) \,\geq\,
\I_{\wt}^{[G]}(0^+) \,,$$
which proves the previous statement. The last thing to prove is that
if $F$ is a distribution function such that $\inf \{ x\,|\, F(x)>0
\} =0$, then
\begin{equation}
\label{eqfini3}
F(0) \,=\, \PP (t(e) = 0)=0 \quad \Longrightarrow \quad \forall \wt\,,\,\,
\I_{\wt}(0^+)\,=\, +\infty\,.
\end{equation}
We consider such a distribution function $F$. We want to compare $F$ with a
Bernoulli distribution of parameter $p$ very close to $1$. For a fixed $p$
(as close to $1$ as we will need), there exists $\eta(p)>0$ such that
$F(\eta(p))<1-p$, because $F(0)=0$ and $F$ is right continuous. We denote by $(t^{[p]}(e))$ the i.i.d. family of Bernoulli
variables of parameter $p$ indexed by the edges, and by $\tau^{[p]}$ the
maximal flow corresponding to these capacities. Then
$$\tau(n\widetilde{A}, h(n))  \,\geq\, \eta(p)\, \tau^{[p]}(n\widetilde{A}, h(n))
\,.  $$
It is proved in section 3 of \cite{Theret:small} (the proof is written for
a straight cylinder $\wt = 0$ and for the variable $\phi$, but it can be
directly adapted to a tilted box and the variable $\tau$ - notice that the
factor $h(n)$ disappears) that there exists a constant $c$ such that for all
$\gamma>0$ 
\begin{align*}
\PP \left( \frac{\tau^{[p]}(n\widetilde{A}, h(n))}{nl(\widetilde{A})} \leq \frac{1}{2}
\right) & \,\leq\, \exp \left(-nl(\widetilde{A}) \left[ \frac{\gamma}{2}
    -\log c - \log (p + (1-p) e^\gamma) \right]  \right)\,.
\end{align*}
Thus for any fixed $R$ (very large, thus $\log c <R$), we can choose
$\gamma= 6 R$, and then $p(R)$ close enough to $1$ to obtain that  $\log (p
+ (1-p) e^{6R})<R$, thus
$$\PP \left( \frac{\tau^{[p(R)]}(n\widetilde{A}, h(n))}{nl(\widetilde{A})}
  \leq \frac{1}{2} \right) \,\leq\, e^{-Rnl(\widetilde{A})}  \,.$$
Finally, for any fixed $R$, for a fixed $\eps$ small enough to have $\eps /
\eta(p(R)) \leq 1/2$, we obtain
\begin{align*}
\PP \left(  \frac{\tau(n\widetilde{A}, h(n))}{nl(\widetilde{A})} \leq \eps
\right) & \,\leq\, \PP \left( \frac{\tau^{[p(R)]}(n\widetilde{A},
    h(n))}{nl(\widetilde{A})} \leq \frac{\eps}{\eta(p(R))} \ \right)\\
& \,\leq\,  e^{-Rnl(\widetilde{A})}\,,
\end{align*}
thus $\I_{\wt} (\eps) \geq R$ for such small $\eps$, which implies that
$\I_{\wt} (0^+) \geq R $ for all $R$. This ends the proof of equation
(\ref{eqfini3}), and thus the proof of Lemma \ref{lemfini}.
\end{dem}
We come back to the proof of Lemma \ref{continuiteK}. We recall that $\wK
(\lambda) = \inf \{ g_{\lambda} (\wt) \,|\, \wt \in \ad(\overline{\D})
\}$. Since $g_{\lambda}$ is l.s.c. and $ad(\overline{\D})$ is compact, there
exists $\wt_{\lambda} \in  \ad(\overline{\D})$ (maybe not unique) such that
$\wK (\lambda) \,=\, g_{\lambda} (\wt_{\lambda})$. If
$$ \lambda \cos (\wt_{\lambda} - \theta)  < \delta  (|\cos \wt| + |\sin \wt|)\,,  $$
then
$$ \wK (\lambda) \,=\, \I_{\wt_{\lambda}} ( \lambda \cos (\wt_{\lambda} -
\theta)  ^+) \,=\, + \infty \,=\, \wK(\lambda ^-) \,. $$
If
$$ \lambda \cos (\wt_{\lambda} - \theta)  > \delta  (|\cos \wt| + |\sin
\wt|)\,,  $$
then
\begin{align*}
\wK(\lambda ^-) &\,=\, \lim_{\eps \rightarrow 0,\, \eps>0} \wK(\lambda -
\eps) \\
& \,\leq\,\lim_{\eps \rightarrow 0,\, \eps>0} \frac{1}{\cos(\wt_{\lambda} -\theta)}\I_{\wt_{\lambda}}
((\lambda - \eps) \cos (\wt_{\lambda} - \theta)  ^+  )\\
& \,\leq\,   \frac{1}{\cos(\wt_{\lambda} -\theta)}\I_{\wt_{\lambda}}
(\lambda \cos (\wt_{\lambda} - \theta)  ^-  )\\
&\,\leq\,  \frac{1}{\cos(\wt_{\lambda} -\theta)}\I_{\wt_{\lambda}}
(\lambda \cos (\wt_{\lambda} - \theta)  ^+  ) \,=\, \wK(\lambda) \,,
\end{align*}
since $\I_{\wt}$ is continuous on $]\delta  (|\cos \wt| + |\sin
\wt|), +\infty[$, thus $\wK(\lambda ^-) = \wK(\lambda)$. We suppose that
$$ \lambda \cos (\wt_{\lambda} - \theta)  = \delta  (|\cos \wt_{\lambda}| + |\sin
\wt|_{\lambda})\,,  $$
which is the only non-obvious case, and thus
$$  \wK (\lambda) \,=\,  \frac{1}{\cos(\wt_{\lambda} -\theta)}\I_{\wt_{\lambda}}
(\delta  (|\cos \wt_{\lambda}| + |\sin \wt_{\lambda}|)   ^+  )\,. $$
If $\PP(t(e)=\delta) =0$, by Lemma \ref{lemfini} we know that the previous
quantity is infinite, thus $\wK(\lambda^-)= +\infty = \wK(\lambda)$. For
the rest of the proof, we suppose that $\PP(t(e) = \delta) >0$. Still by
Lemma \ref{lemfini}, we know in this case that for all $\wt$, $\I_{\wt}
(\delta  (|\cos \wt| + |\sin \wt|)^+ ) <+\infty$. If $\lambda =
\delta_{\theta, h}$, we have nothing to prove, thus we suppose that
$\lambda > \delta_{\theta,h}$ (it implies that $\lambda >0$). We suppose that the following
property $\mathcal{P}$ holds: there exists a sequence $(\wt_n)_{n\in \NN}$
such that
$$ \mathcal{P} \left\{ \begin{array}{l} (i) \,\, \lim_{n\rightarrow \infty}\wt_n \,=\,
  \wt_{\lambda}\,,\\ (ii)\,\, \forall n\in \NN\,,\,\, \delta \frac{|\cos
    \wt_n| + |\sin \wt_n|}{\cos(\wt_n - \theta)} \,<\, \lambda \,=\, \delta \frac{|\cos \wt_{\lambda}| + |\sin
    \wt_{\lambda}|}{\cos(\wt_{\lambda} - \theta)}  \,, \\ (iii)\,\, \limsup_{n\rightarrow
    \infty} \I_{\wt_n} (\delta (|\cos \wt_n| + |\sin
    \wt_n|)^+ ) \,\leq\, \I_{\wt_{\lambda}} (\delta (|\cos \wt_{\lambda}| + |\sin
    \wt_{\lambda}|) ^+)\,.  \end{array}\right. $$
We consider a given $\eta>0$. For $n_0$ large enough, we have
$$ \frac{1}{\cos(\wt_{n_0} -\theta)}\I_{\wt_{n_0}} (\delta (|\cos
\wt_{n_0}| + |\sin \wt_{n_0}|) ^+  ) \,\leq\,  \frac{1}{\cos(\wt_{\lambda} -\theta)}\I_{\wt_{\lambda}} (\delta (|\cos
\wt_{\lambda}| + |\sin \wt_{\lambda}|) ^+  ) + \eta\,.  $$
Moreover, there exists $\eps_0 >0$ such that
$$ \delta \frac{|\cos \wt_{n_0}| + |\sin
    \wt_{n_0}|}{\cos(\wt_{n_0} - \theta)} \,\leq\, \lambda -\eps_0\,, $$
and for all $\eps\leq\eps_0$, since $\I_{\wt_{n_0}}$ is non increasing, we
obtain that
\begin{align*}
\wK (\lambda -\eps) & \,\leq\, \frac{1}{\cos(\wt_{n_0}
  -\theta)}\I_{\wt_{n_0}} ((\lambda - \eps) \cos(\wt_{n_0} - \theta)^+ ) \\
& \,\leq\,  \frac{1}{\cos(\wt_{n_0} -\theta)}\I_{\wt_{n_0}} (\delta (|\cos
\wt_{n_0}| + |\sin \wt_{n_0}|) ^+  )\\
& \,\leq\,  \frac{1}{\cos(\wt_{\lambda} -\theta)}\I_{\wt_{\lambda}} (\delta (|\cos
\wt_{\lambda}| + |\sin \wt_{\lambda}|) ^+  ) + \eta \,=\, \wK(\lambda) + \eta\,.
\end{align*}
We conclude that $\wK(\lambda^-) \leq \wK(\lambda)$, so $\wK(\lambda^-) =
\wK(\lambda)$, and this ends the proof of Lemma \ref{continuiteK}.

The last thing we have to do is to prove the property
$\mathcal{P}$. Obviously, property $(ii)$ is linked with the monotonicity of
the function
$$\Gamma : \wt \mapsto \frac{|\cos \wt | + |\sin \wt |}{\cos(\wt - \theta)}
\,. $$
If $\theta \in \{k\pi/4\,|\, k\in \NN \} $, we will prove in the next
paragraph, see Lemma \ref{chapitre6Kcasdroit}, that $\wK (\lambda)=\I_{\theta}
(\lambda ^+)$; since it is obvious that in this case $\delta_{\theta,h} =
\delta$, the continuity of $\wK$ possibly except at $\delta_{\theta,h}$ is
already known. We suppose that $\theta \notin \{k\pi/4\,|\, k\in \NN \} $,
and by symmetry we can suppose that $\theta \in ]0,\pi/2[ \smallsetminus \{
\pi/4 \}$. It is obvious
(see the factor $\cos (\wt-\theta) ^{-1}$)
that for all $\wt \in [\theta-\pi/2, \theta +\pi/2] \smallsetminus
[0,\pi/2]$, we have $\Gamma(\wt) > \inf_{[0,\pi/2]} \Gamma$, so $\argmin
\Gamma \in [0,\pi/2]$. Similarly,
$\wt_{\lambda} \in [0,\pi/2] $ too. For all $\wt \in [0,\pi/2]$, we can write
$$ \Gamma (\wt)  \,=\, \frac{\cos \wt + \sin \wt}{ \cos (\wt - \theta)}
\,=\,  \frac{1+ \tan \wt}{ \cos \theta + \sin \theta \tan \wt}\,. $$
We deduce from this equality that $\Gamma$ is strictly monotone on
$[0,\pi/2]$: strictly increasing (resp. decreasing) if $\theta \in
]0,\pi/4[$ (resp. $\theta \in ]\pi/4, \pi/2[$ ), and thus
$\argmin \Gamma = 0$ (resp. $\argmin \Gamma =\pi/2$). We consider the
case $\theta \in ]0,\pi/4[$, the study of the case $\theta \in
]\pi/4,\pi/2[$ being similar. We know that $\wt_{\lambda} \in ]0,\pi/2]$,
because $\wt_{\lambda} =0$ implies that $\lambda = \delta_{\theta, h}$,
and we excluded this case. Thus we can consider a strictly increasing
sequence $(\wt_n)_{n\in\NN}$ such that $\wt_n \in ]0,\pi/2[$ for all $n$
and $(i):\, \lim_{n\rightarrow \infty} \wt_n =\wt_{\lambda}$ is
satisfied. Since $\Gamma$ is strictly increasing on $[0,\pi/2]$, we know
that such a strictly increasing sequence $(\wt_n)_{n\in \NN}$ satisfies the
hypothesis $(ii)$. To prove that $(iii)$ also holds, we need Lemmas
\ref{lemconvexiteLambda} and \ref{lemfini}. The function
$$ \Lambda_{\delta} : \vec{v} \mapsto \Lambda (\delta, \vec{v}) $$
is finite on $(\RR^+)^2\smallsetminus \{(0,0)\}$ on the hypothesis
$F(\delta)>0$ we did (see Lemma \ref{lemfini}), and it is convex (see
Lemma \ref{lemconvexiteLambda}), so it is continuous on the interior of
$(\RR^+)^2\smallsetminus \{(0,0)\}$. If $\wt_{\lambda} \neq \pi/2$, it
proves $(iii)$. We suppose $\wt_{\lambda} =\pi/2$. Let $\vec{u} = (0,1)$,
$\vec{v}_n = (1/\tan\wt_n, 0) $ and $\vec{w}_n = (1/\tan\wt_n, 1) =
\vec{u} + \vec{v}_n$. By equation (\ref{eqconvexite}) for $\lambda =
\delta$ we have for all $n\in \NN$
$$ \frac{1}{\sin \wt_n} \I_{\wt_n} (\delta (|\cos \wt_n| + |\sin
\wt_n|)^+) \,\leq\, \I_{\pi/2} (\delta^+) +\frac{1}{\tan\wt_n}
\I_{0}(\delta^+)  \,,$$ 
and sending $n$ to infinity we exactly obtain $(iii)$, so the property
$\mathcal{P}$ is proved.
\end{dem}


We prove finally the property stated in Remark \ref{chapitre6pgdcasdroit},
in fact a property a little bit more general:
\begin{lem}
\label{chapitre6Kcasdroit}
If $\theta \in \{k\pi/4\,|\, k\in \NN \}$, then
$$ \K \,=\, \J_{\theta}  \,.$$
\end{lem}

\begin{dem}
We fix a $\theta \in \{k\pi/4\,|\, k\in \NN\}$. We know that $\nu_{\theta}
= \eta_{\theta,h}$ for such a $\theta$ (see Remark 2.11 in
\cite{RossignolTheret09}), so it is sufficient to prove that
$$ \forall \lambda \geq 0 \,,\qquad \wK (\lambda) \,=\,
\I_{\theta} (\lambda ^+)\,. $$
Since $\theta\in \ad(\overline{\D})$, it is equivalent to prove that
$$ \forall \lambda \geq 0\,,\,\, \forall \wt \in [\theta -\pi/2,
\theta+\pi/2]\,, \qquad 
\I_{\theta}(\lambda ^+ )\,\leq\,\frac{1}{\cos (\wt - \theta)}
\I_{\wt}(\lambda \cos (\wt - \theta)^+)\,.  $$
Let $\wt \in ]\theta -\pi/2, \theta + \pi/2 [$. We use the same notations
as in Lemma \ref{chapitre6trigI}. We consider the non
degenerate triangle $(abc)$ such that $\wt_c = \theta + \pi $ (so
$\cyl_c(n)$ is a straight cylinder in the case $\theta =0$), $\wt_b = \max( \wt
, 2\theta - \wt)$, $\wt_a
=\min( \wt , 2\theta - \wt) $, $l(ab)=1$ and $l(ac)=l(bc)= (2 \cos(\wt -
\theta))^{-1}$. Since the graph is invariant by a symmetry of axis $((0,0),
(\cos \theta, \sin \theta))$ (respectively $((0,0), (1,1))$), we know that
$\I_{\wt_a}= \I_{\wt_b}$ (respectively $\I_{\wt_c} = \I_{\theta}$). Then
Lemma \ref{chapitre6trigI} applied with $\alpha = 1/2$ states that for all $\lambda
\geq 0$,
$$ \I_{\theta} \left( \lambda ^+ \right) \,\leq\, \frac{1}{\cos (\wt -
  \theta)} \I_{\wt} (\lambda \cos(\wt - \theta) ^+) \,. $$
The inequality remains obviously valid for $\wt \in \{\theta + \pi/2,
\theta -\pi/2\}$, since we have seen in Remark \ref{chapitre6reminfini} that the
right hand side of the previous inequality equals $+\infty$ in this
case. This ends the proof of Lemma \ref{chapitre6Kcasdroit}.
\end{dem}


\subsection{Lower bound}

We have to prove that for all open subset $\mathcal{O}$ of $\RR^+$, we have
$$ \liminf_{n\rightarrow \infty} \frac{1}{nl(A)} \log \PP
\left[\frac{\phi_n}{nl(A)}  \in \mathcal{O} \right]\, \geq\,  -\inf_{\lambda \in
  \mathcal{O}} \K(\lambda) \,.$$
Classically, it suffices to prove the local lower bound:
\begin{equation}
\forall a\in \RR^+\,,\,\, \forall \varepsilon >0 \qquad
\liminf_{n\rightarrow \infty} \frac{1}{nl(A)} \log \PP
\left[\frac{\phi_n}{nl(A)}  \in [a-\varepsilon, a+\varepsilon] \right]
\,\geq\,  -\K(a)\,.
\end{equation}
If $\K(a) = +\infty$, the result is obvious, so we suppose that
$\K(a)<+\infty$. For all $\eta < \varepsilon$, we have
\begin{align}
\liminf_{n\rightarrow \infty} \frac{1}{nl(A)} \log &\PP
\left[\frac{\phi_n}{nl(A)}  \in [a-\varepsilon, a+\varepsilon] \right]
\nonumber \\
&\,\geq\, \liminf_{n\rightarrow \infty} \frac{1}{nl(A)} \log \left( \PP
\left[\frac{\phi_n}{nl(A)} \leq a+\eta \right] -\PP \left[
  \frac{\phi_n}{nl(A)} \leq a-\varepsilon \right] \right)\,. \label{chapitre6som}
\end{align}
From the strict decreasing of $\K$ (see Lemma \ref{chapitre6strict}), we deduce
that for all $a\in \RR^+$ such that $\K(a)<\infty$, for all positive $\eta$ and
$\varepsilon$, we have
\begin{equation}
\label{chapitre6strictbis}
\inf_{\wt \in \D} \frac{1}{\cos (\wt - \theta)} \I_{\wt} ((a+\eta)
\cos(\wt - \theta)^-) \,<\, \inf_{\wt \in \D} \frac{1}{\cos (\wt - \theta)}
\I_{\wt} ((a-\varepsilon) \cos(\wt - \theta)^+) \,. 
\end{equation}
Indeed, for all positive $\eta$, we have
$$ \inf_{\wt \in \D} \frac{1}{\cos (\wt - \theta)} \I_{\wt} ((a+\eta)
\cos(\wt - \theta)^-) \,\leq\, \K(a) \,<\, \K(a-\varepsilon) \,.$$
Then thanks to (\ref{chapitre6eq1'}), (\ref{chapitre6eq2}) and (\ref{chapitre6strictbis}), we know
that the second term in the sum appearing in (\ref{chapitre6som}) is negligible
compared to the first one, so we obtain that
\begin{align*}
\liminf_{n\rightarrow \infty} \frac{1}{nl(A)} \log & \PP
\left[\frac{\phi_n}{nl(A)}  \in [a-\varepsilon, a+\varepsilon] \right]\\
& \,\geq\, -\inf_{\wt \in \D}  \frac{1}{\cos (\wt - \theta)} \I_{\wt} ((a+\eta)
\cos(\wt - \theta)^-)\\
& \,\geq\,  -\inf_{\wt \in \D} \lim_{\varepsilon'\rightarrow 0}
\frac{1}{\cos (\wt - \theta)} \I_{\wt} ((a+\eta) \cos(\wt - \theta) -
\varepsilon') \,.
\end{align*} 
Sending $\eta$ to zero, we obtain that
\begin{align*}
\liminf_{n\rightarrow \infty} \frac{1}{nl(A)} \log &\PP
\left[\frac{\phi_n}{nl(A)}  \in [a-\varepsilon, a+\varepsilon] \right]\\
& \,\geq\, -\liminf_{\eta \rightarrow 0} \inf_{\wt \in \D}
\lim_{\varepsilon'\rightarrow 0} \frac{1}{\cos (\wt - \theta)} \I_{\wt}
((a+\eta) \cos(\wt - \theta) - \varepsilon')\\
& \,\geq\,-  \inf_{\wt \in \D}\lim_{\eta \rightarrow 0}
\lim_{\varepsilon'\rightarrow 0} \frac{1}{\cos (\wt - \theta)} \I_{\wt}
((a+\eta) \cos(\wt - \theta) - \varepsilon') \\
& \,\geq\,-  \inf_{\wt \in \D} \frac{1}{\cos (\wt - \theta)} \I_{\wt}
(a \cos(\wt - \theta) ^+) \,,
\end{align*}
and so the local lower bound is proved.





\subsection{Upper bound}

\subsubsection{Upper large deviations}
\label{seculd}

To handle the upper large deviations, we shall use the following result:
\begin{lem}
\label{chapitre6up}
Let $A$ be a non-empty line-segment in $\RR^2$, with Euclidean length
$l(A)$. Let $\theta\in[0,\pi[$ be such that $(\cos\theta,\sin\theta)$
is orthogonal to the hyperplane spanned by $A$ and $h:\NN\rightarrow
\RR^+$. We suppose that $\mathbf{(H1)}$, $\mathbf{(H2)}$, $\mathbf{(F2)}$,
$\mathbf{(FH1)}$ and either $\mathbf{(F4)}$ or  $\mathbf{(H3)}$ hold. Then
for all $\lambda > \eta_{\theta,h}$ we have
\begin{equation}
\label{eqlemsup}
 \limsup_{n\rightarrow \infty} \frac{1}{nl(A)} \log \PP \left[
  \frac{\phi_n}{nl(A)} \geq \lambda  \right] \,= \,-\infty \,. 
\end{equation}
\end{lem}

In fact, we have a stronger result, if $F$ admits an exponential moment:
\begin{lem}
\label{chapitre6thmupper}
Let $A$ be a non-empty line-segment in $\RR^2$, with Euclidean length
$l(A)$. Let $\theta\in[0,\pi[$ be such that $(\cos\theta,\sin\theta)$
is orthogonal to the hyperplane spanned by $A$ and $h:\NN\rightarrow
\RR^+$. We suppose that $\mathbf{(H1)}$, $\mathbf{(H2)}$, $\mathbf{(F2)}$,
$\mathbf{(FH1)}$ and $\mathbf{(F4)}$ hold. Then for all $\lambda
>\eta_{\theta,h}$, we have
$$ \liminf_{n\rightarrow \infty} \frac{-1}{nl(A) h(n)} \log \PP \left[
  \phi(nA,h(n)) \geq \lambda n l(A) \right] \,>\, 0 \,. $$
The upper large deviations are thus of volume order.
\end{lem}

Obviously, Lemma \ref{chapitre6thmupper} implies Lemma \ref{chapitre6up} in
the case where the condition $\mathbf{(F4)}$ is satisfied, since
$\lim_{n\rightarrow \infty} h(n) =+\infty$ by $\mathbf{(H1)}$. We do not
present a complete version of the proof of Lemma \ref{chapitre6thmupper}:
it is simply a modification of the proofs of Theorem 2 in
\cite{TheretUpper} (the part concerning the positivity of the rate
function, section 3.7) and Theorem 4 in \cite{Theret:uppertau}, and it can
be found in Part 3, Chapter 6, section 5 of \cite{Theret:these}. The common
idea of these proofs is the following. We consider the cylinder
$\cyl(NA,h(N))$, and divide it into slabs of height $2h(n)$, i.e.,
translates of $\cyl(NA,h(n))$, for $n$ a lot smaller than $N$. If
$\phi(NA,h(N))$ is big, it implies that the maximal flow from the top to
the bottom of each slab is big too, and we have of order $h(N)$ such slabs
for a fixed $n$. It implies roughly that
\begin{equation}
\label{eqsup}
\PP\left[ \phi(NA,h(N)) \geq \lambda Nl(A) \right] \,"\leq"\, \PP
\left[ \phi(NA,h(n)) \geq \lambda Nl(A)  \right]^{p h(N)}\,,
\end{equation}
for some constant $p$. We divide then each slab
into disjoint translates of $\cyl(nA,h(n))$, and we can compare the maximal
flow from the top to the bottom of the slab with the sum of the variables
$\tau^i_n$ into these small cylinders. Roughly speaking, we obtain that
$$ \PP \left[ \phi(NA,h(n)) \geq \lambda Nl(A)  \right] \,"\leq" \, \PP
\left[ \sum_i \tau^i_n \geq \lambda Nl(A)  \right] \,,$$
Under the hypothesis $\mathbf{(F4)}$, Cram\'er's Theorem in $\RR$ states
that $\PP [\sum_i \tau^i_n \geq \lambda Nl(A)]$ decays exponentially fast
with $N$ for any $\lambda > \nu_{\theta} = \lim_{n\rightarrow \infty}
\tau^i_n / (nl(A))$, thus we obtain that
$$ \PP \left[ \phi(NA,h(N)) \geq \lambda Nl(A)  \right] \,\leq \, p'
e^{- p'' h(N)Nl(A)}\,, $$
for other constants $p'$ and $p''$. The only adaptation we have to do is to
take into account the fact that under hypothesis $\mathbf{(FH1)}$, the limit
$\eta_{\theta, h}$ of $\phi(NA,h(N))/(Nl(A))$ is equal to
$\nu_{\wt_0}/\cos(\wt_0 - \theta)$ for some $\wt_0$. Thus we divide
$\cyl(NA,h(N))$ into slabs orthogonal to $\vec{v}(\wt_0)$ instead of slabs
orthogonal to $\vec{v}(\theta)$. Thus we compare $\phi(NA,h(N))$ with $h(N)$ sums of $N l(A) /\cos(\wt_0 - \theta)$
terms equal in law with $\tau_n(\wt_0)$, the maximal flow from the upper
half part to the lower half part of the boundary of a box of size $n\times
h(n)$ oriented towards the direction $\wt_0$. We conclude again thanks to
 Cram\'er's Theorem in $\RR$.

If $\mathbf{(F4)}$ is not satisfied, we cannot use Cramer's
Theorem. However, we can perform the division of $\cyl(NA,h(N))$ into slabs
orthogonal to $\wt_0$, and thus obtain an equation very close to
(\ref{eqsup}):
\begin{equation}
\label{eqhgrand}
 \PP\left[ \phi(NA,h(N)) \geq \lambda Nl(A) \right] \,"\leq"\, \PP
\left[ \tau^i_{N,n} \geq \lambda Nl(A)  \right]^{p h(N)}\,,
\end{equation}
where $\tau^i_{N,n}$ is the maximal flow from the upper half part to the
lower half part of the boundary of a slab. Thus, if $\mathbf{(H3)}$ holds
instead of $\mathbf{(F4)}$, equation (\ref{eqhgrand}) leads to the
conclusion of Lemma \ref{chapitre6up}.

\begin{rem}
\label{remcondsup}
The hypotheses $\mathbf{(H3)}$ or $\mathbf{(F4)}$ may not be optimal, but a
simple example shows why we need such kind of hypotheses. We consider that
the capacity of an edge is distributed according to the Pareto law of
parameters $p$ and $1$, i.e., the probability that an edge has a capacity bigger
than $t \geq 1$ is equal to $t^{-p}$. We consider the rectangle $A=[0,1]\times
\{0\}$, and the maximal flow $\phi(nA,h(n))$ from the top to the bottom of
the cylinder $[0,n]\times [-h(n),h(n)]$. If all the vertical edges $(e_i,
i=1,...,2h(n))$ (we suppose $h(n)\in \NN$ for simplicity) in the box that are included in the segment $\{1\} \times [-h(n),h(n)]$ have a
capacity bigger than $\lambda n$ for a fixed $\lambda$, then
$\phi(nA,h(n))$ is bigger than $\lambda n$. We obtain:
$$ \PP \left[ \phi(nA,h(n)) \geq \lambda n \right] \,\geq\, \PP\left[
  \forall i=1,...,2h(n)\,,\, t(e_i)\geq \lambda n \right] \,\geq\, (\lambda n)^{-2ph(n)}
\,.$$
If $ h(n) \log n $ is not large compared to $n$, in the sense that $h(n)
\log n / n$ does not converge towards $+\infty$, then equation
(\ref{eqlemsup}) is not satisfied.
\end{rem}

\subsubsection{End of the proof of Theorem \ref{thmPGDphi}}

For this last section, we impose $\mathbf{(H1)}$, $\mathbf{(H2)}$,
$\mathbf{(F1)}$, $\mathbf{(F2)}$, $\mathbf{(FH1)}$, $\mathbf{(FH2)}$ and
either  $\mathbf{(F4)}$ or  $\mathbf{(H3)}$. Let $\F$ be a closed subset of $\RR^+$. We want to prove that
$$ \limsup_{n\rightarrow \infty} \frac{1}{nl(A)} \log \PP \left[
  \frac{\phi_n}{ nl(A)} \in \F \right] \,\leq\, -\inf_{\lambda \in \F}
\K(\lambda) \,.  $$
If $\eta_{\theta,h}$ belongs to $\F$, then according to Corollary
\ref{chapitre6corollaire}, we know that
$$ \lim_{n\rightarrow \infty}\PP \left[ \frac{\phi_n}{nl(A)} \in \F
\right] \,=\, 1 \,,$$
and so
$$\limsup_{n\rightarrow \infty} \frac{1}{nl(A)} \log \PP \left[
  \frac{\phi_n}{ nl(A)} \in \F \right] \,=\, 0 \,=\,  -\inf_{\lambda \in \F}
\K(\lambda)\,, $$
because $\K$ is non-negative, and $\K(\eta_{\theta, h}) =0$.
Let us suppose that $\eta_{\theta, h}$ does not belong to $\F$. The following
proof is similar to the one of the upper bound in \cite{RossignolTheret08b}. We define $f_1 = \sup (\F\cap [0, \eta_{\theta, h}])$ and $f_2 =(
\inf \F \cap [\eta_{\theta,h}, +\infty[)$. We suppose here that $\F\cap [0,
\eta_{\theta, h}]$ and $\F \cap [\eta_{\theta,h}, +\infty[$ are non empty,
because it is the most complicated case (if one of these two sets is
empty, part of the following study is sufficient). Since $\F$ is closed,
we know that $f_1 < \eta_{\theta,h}$ and $f_2> \eta_{\theta, h}$. Then
\begin{align*}
\limsup_{n\rightarrow \infty} \frac{1}{nl(A)} \log & \PP \left[
  \frac{\phi_n}{ nl(A)} \in \F \right]\\ & \,\leq\, \limsup_{n\rightarrow
  \infty} \frac{1}{nl(A)} \log \left( \PP \left[ \frac{\phi_n}{
      nl(A)} \leq f_1 \right] + \PP\left[  \frac{\phi_n}{
      nl(A)} \geq f_2 \right] \right)\,.
\end{align*}
On one hand, by (\ref{chapitre6eq2}), we know that
$$ \limsup_{n\rightarrow \infty }
\frac{1}{nl(A)} \log \PP [\phi_n \leq f_1 n l(A)] \,\leq\, -\K(f_1)
\,. $$
On the other hand, if we refer to Lemma \ref{chapitre6up}, we know that
$$ \limsup_{n\rightarrow \infty } \frac{1}{nl(A)} \log \PP [\phi_n \geq
f_2 n l(A)] \,=\, - \infty \,. $$ 
If $\K(f_1) = + \infty$, we have
$$ \limsup_{n\rightarrow \infty} \frac{1}{nl(A)} \log \PP \left[
  \frac{\phi_n}{ nl(A)} \in \F \right] \,=\, - \infty \,=\, - \inf_\F \K
\,,$$
because $\K$ is infinite on $[0,f_1]$ ($\K'$ is non-increasing) and on
$[f_2, +\infty[$, so on $\F$. If $\K(f_1)<\infty$, we have
$$  \limsup_{n\rightarrow \infty} \frac{1}{nl(A)} \log \PP \left[
  \frac{\phi_n}{ nl(A)} \in \F \right] \,\leq\, -\K(f_1) \,=\, - \inf_\F \K
\,,  $$
because $\K$ is non-increasing on $[0, f_1]$ and infinite on $[f_2,
+\infty[$. So the upper bound is proved.



\def\cprime{$'$}

\end{document}